\documentclass[final]{siamltex}

\usepackage{amssymb,amsmath, amsfonts, amsbsy, amsbsy, latexsym, color, epsfig, graphicx}

\newcommand{\bitem}{\begin{itemize}}
\newcommand{\eitem}{\end{itemize}}
\newcommand{\benum}{\begin{enumerate}}
\newcommand{\eenum}{\end{enumerate}}
\newcommand{\beq}{\begin{equation}}
\newcommand{\eeq}{\end{equation}}

\newcommand{\Sp}{\mbox{supp}}
\newcommand{\rd}{\mbox{round}}
\newcommand{\sgn}{\mbox{sgn}}

\def\NN{\mathbb{N}}
\def\ZZ{\mathbb{Z}}
\def\RR{\mathbb{R}}
\def\ZZ{\mathbb{Z}}

\newcommand{\bZ}{{\mathbb Z}}
\newcommand{\bR}{{\mathbb R}}
\newcommand{\bN}{{\mathbb N}}

\def\cC{{\mathcal{C}}}

\def\cR{{\mathcal{R}}}

\def\cI{{\mathcal{I}}}
\def\cL{{\mathcal{L}}}
\def\cAS{{\mathcal{A}\mathcal{S}}}

\newcommand{\iu}{i}         

\newcommand{\df}{\mathsf{d}}
\newcommand{\dm}{\mathsf{M}}        

\newcommand{\dn}{\mathsf{N}}        

\newcommand{\dmcg}[1]{\Gamma_{#1}}  
\newcommand{\dmfcg}[1]{\Omega_{#1}} 

\newcommand{\C}{\mathbb{C}}    
\newcommand{\Z}{\mathbb{Z}}    

\newcommand{\dR}{\mathbb{R}^d}

\newcommand{\dZ}{\mathbb{Z}^d}

\newcommand{\dHH}[1]{{H^{#1}(\mathbb{R}^d)}}         

\newcommand{\llll}{l}

\newcommand{\vk}{\mathsf{k}}

\newcommand{\vn}{\mathsf{n}}

\newcommand{\ta}{a}
\newcommand{\tb}{b}
\newcommand{\tu}{u}
\newcommand{\tv}{v}
\newcommand{\tw}{w}

\newcommand{\mphi}{s}
\newcommand{\mpsi}{r}
\newcommand{\mi}{\ell}  

\newcommand{\sd}{{\mathrm{S}}}  
\newcommand{\tz}{{\mathrm{T}}}  

\title{A Unitary Extension Principle for Shearlet Systems}

\author{Bin Han\thanks{Department of Mathematical and Statistical Sciences, University of Alberta,
Edmonton, Alberta, Canada T6G 2G1 ({\tt bhan@math.ualberta.ca}).
B.H. is supported in part by NSERC Canada.}
\and Gitta Kutyniok\thanks{Institute of Mathematics, University of Osnabr\"uck, 49069 Osnabr\"uck, Germany
({\tt kutyniok@math.uni-osnabrueck.de}). G.K.  would like to thank the Department of Mathematics at the National
University of Singapore for its hospitality and support during her
visit which enabled completion of a significant part of this paper.}
\and Zouwei Shen\thanks{Department of Mathematics, National University of Singapore, Singapore 117543
({\tt matzuows@nus.edu.sg}). Z.S. is supported in part by Grant
R-146-000-113-112 at the National University of Singapore.}}

\begin{document}

\maketitle

\begin{abstract}
In this paper, we first introduce the concept of an adaptive MRA (AMRA) structure which is
a variant of the classical MRA structure suited to the main goal of a fast flexible
decomposition strategy adapted to the data at each decomposition level. We then
study this novel methodology for the general case of affine-like systems, and derive a
Unitary Extension Principle (UEP) for filter design. Finally, we apply our results to the
directional representation system of shearlets. This leads to a comprehensive theory for
fast decomposition algorithms associated with shearlet systems which encompasses tight
shearlet frames with spatially compactly supported generators within
such an AMRA structure. Also shearlet-like systems associated with parabolic scaling and
unimodular matrices optimally close to rotation as well as 3D shearlet systems are studied
within this framework.
\end{abstract}

\begin{keywords}
Affine systems, fast decomposition algorithm, shearlets, multiresolution analysis, tight frames
\end{keywords}

\begin{AMS}
Primary 42C40; Secondary 42C15, 65T60, 65T99, 94A08
\end{AMS}

\pagestyle{myheadings}
\thispagestyle{plain}
\markboth{B. HAN, G. KUTYNIOK, AND Z. SHEN}{A UNITARY EXTENSION PRINCIPLE FOR SHEARLET SYSTEMS}

\section{Introduction}

Wavelets are nowadays indispensable as a multiscale encoding system
for a wide range of more theoretically to more practically oriented
tasks, since they provide optimal approximation rates for smooth
1-dimensional data exhibiting singularities. The facts that they
provide a unified treatment in both the continuous as well as
digital setting and that the digital setting admits a
multiresolution analysis leading to a fast spatial domain decomposition
were essential for their success. It can however be shown that
wavelets -- although perfectly suited for isotropic structures -- do
not perform equally well when dealing with anisotropic phenomena.

This fact has motivated the development of various types of
directional representation systems for 2-dimensional data that are
capable of resolving edge- or curve-like features which precisely
separate smooth regions in sparse way, of which examples are
contourlets \cite{DV05} and curvelets \cite{CD04,CD05a,CD05b}. All
these multiscale variants offer different advantages and
disadvantages, however, neither of them provides a unified treatment
of the continuous and digital setting. Curvelets, for instance, are
known to yield tight frames but the digital curvelet transform is
not designed within the curvelet-framework and hence, in particular,
is not covered by the available theory \cite{CDDY06}.

About three years ago, a novel representation system -- so-called shearlets
-- has been proposed \cite{GKL06,KL09}, which possesses the same
favorable approximation and sparsity properties as the other
candidates (see \cite{GL07,DK08b,KL10}) of whom curvelets are perhaps the
most advanced ones. One main point in comparison with curvelets is
the fact that angles are replaced by slopes when parameterizing
directions which significantly supports the treating of the digital
setting. A second main point is that shearlets fit within the
general framework of affine-like systems, which provides an
extensive mathematical machinery. Thirdly, it is shown in
\cite{DKS09b} that shearlets  -- in addition to the
aforementioned favorable properties -- provide a unified treatment
for the continuous and digital world similar to wavelets.

Recently, several researchers \cite{GLLWW06,KS08,Lim09} have provided
approaches to introduce an MRA structure with accompanying fast
spatial domain decomposition for shearlets. This would
establish shearlets as the directional representation system which
provides the full range of advantageous properties for 2D data which
wavelets provide in 1D. However, in the previous approaches, does
either the MRA structure not lead to a tight frame, infinitely many
filters make an implementation difficult, or the MRA structure is
not faithful to the continuum transform. Approaches to extend the
present shearlet theory to higher dimensions were also already
undertaken, however, for now, only with continuous parameters
\cite{DST09}.

In wavelet theory, the Unitary Extension Principle (UEP) introduced in
\cite{RS97a} has proven to be a highly useful methodology for constructing
tight wavelet frames with an associated MRA structure in order to develop
efficient algorithms, for instance, for
frame based image restorations including image deblurring and blind
image deblurring, image inpainting, image denoising, and image
decomposition (see \cite{CCSS,CCSS:NM:09, CCSS:ACM:09, CCS, blind2, blind,COS,
Cai2009,CS:NM:07,CCSS:SISC:03, CRSS:ACHA:04,CSX:ACHA:07}).
These applications of frame-based image restorations are also one
motivation for our adventures here.

In this paper we aim at providing an MRA structure for tight
shearlet frames -- and in fact even for more general affine-like
systems encompassing different shearlet systems as special cases --
which exhibits all the favorable properties of MRA structures for
wavelets. We also allow the MRA structure to be
more flexible in the sense of adaptivity than ordinarily considered
in the literature. We further prove sufficient conditions for such
an adaptive MRA in terms of a suited Unitary Extension Principle (UEP) along with
a fast spatial domain decomposition as well as approximation properties.
Surprisingly, our general theory also includes shearlet systems for
3D data and provides fast decomposition algorithms for those.


\subsection{List of Desiderata}
\label{subsec:desiderata}

The theoretical framework for a fast spatial domain decomposition
we will develop in the case of general affine-like systems (cf. Subsection \ref{subsec:AMRA} for
a precise definition), and, in particular, for shearlet systems shall satisfy the following
desiderata:

\bitem
\item[(1)] {\em Adaptivity.} The MRA structure shall allow interactive adaption of the
decomposition procedure to different types of data either beforehand
or during the process.
\item[(2)] {\em Computational Feasibility.} The MRA structure shall be suited to
the fact that it aims for a decomposition algorithm, that means, for
instance, that only a finite range of scales is required.
\item[(3)] {\em Fast Spatial Domain Decomposition.} An accompanying fast transform shall exist
which decomposes digital data into high- and low-frequency parts, in
particular, comprising directionality behavior.
\item[(4)] {\em Tight Frame Property.} This property will ensure stability of the decomposition
and the possibility to employ the adjoint for reconstruction.
\item[(5)] {\em Compactly-Supportedness.} The framework shall allow the utilization of compactly
supported generators to ensure spatial domain localization as well as fast computations, in
particular, in the special case of shearlets.
\item[(6)] {\em Arbitrary Dimension.} The framework shall be applicable to affine-like systems
for arbitrary dimension.
\item[(7)] {\em Faithfulness to Continuum Transform.} In the special case of shearlets, the
transform shall be faithful to the continuum transform in the sense of being accompanied by
a continuum shearlet system with discrete parameters.
\eitem


\subsection{An Adaptive Multiresolution Analysis (AMRA) for Affine-Like Systems}
\label{subsec:AMRA}

The framework of a multiresolution analysis is a well-established
methodology in wavelet theory for deriving a decomposition of data
into low- and high-frequency parts associated with a scaling
function and wavelet which leads to a fast spatial domain
decomposition. However, aiming for a decomposition with respect to general {\em affine-like systems},
which we coin systems being a subset of unions of {\em affine systems}
\[
\bigcup_{\dm \in \Lambda \subseteq GL(d,\RR)} \{\psi^\mi_{\dm^j;\vk}\; :\;
j\in \Z, \vk \in \dZ, \mi=1,\ldots,\mpsi\}, \qquad
\psi^1, \ldots, \psi^\mpsi \in L^2(\bR^d),
\]
where
\begin{equation}\label{dilation:shift}
\psi^\mi_{U; \vk}:=|\det U|^{1/2} \psi^\mi(U\cdot-\vk), \qquad
\vk\in \dZ, U\in GL(d, \RR),
\end{equation}
we are forced to reconsider this approach.

In light of desiderata (1)--(3) from the previous subsection and
continuing the initial observations already made by one of the
authors in \cite{Han09}, we claim a new paradigm for MRA manifested
in the following two requirements to an MRA:
\bitem
\item[(R1)] If one aims for a fast decomposition algorithm, it is sufficient to study nonhomogeneous
systems, i.e., considering a fixed finest decomposition level $J$
instead of analyzing the limit $J \to \infty$. In the same way, the
coarsest level might be defined to be the level $0$ without loss of
generality.
\item[(R2)] It is not necessary to specify beforehand which functions or masks generate the low- and which
generate the high-frequency parts, since the key equation (Unitary
Extension Principle) which needs to be satisfied at each
decomposition level
does not distinguish between them.
\eitem
The main advantage of (R1) is to adapt the theoretical continuum considerations to the
implementation requirements, whereas (R2) allows us to change
the composition of the masks adaptively at each decomposition level.
Therefore we coin this new paradigm for MRA an {\em Adaptive
Multiresolution Analysis (AMRA)}.

The decomposition technique will still rely on subdivision and
transition operators which we in the sequel denote as follows: For a
$d \times d$ invertible integer matrix $\dm$ and a finitely
supported sequence $\ta : \ZZ^d \to \RR$, the {\em subdivision
operator} $\sd_{\ta,\dm} : \llll(\ZZ^d) \to \llll(\ZZ^d)$ and the
{\em transition operator} $\tz_{\ta,\dm} : \llll(\ZZ^d) \to
l(\ZZ^d)$ are defined by
\[
[\sd_{\ta,\dm} \tv](\vn) := |\det(\dm)| \sum_{\vk \in \ZZ^d}
\tv(\vk) \ta(\vn-\dm \vk) \quad \mbox{and} \quad
[\tz_{\ta,\dm} \tv](\vn) :=  \sum_{\vk \in \ZZ^d} \tv(\vk)
\overline{\ta(\vk-\dm\vn)}.
\]
Note that $\widehat{\sd_{\ta,\dm} \tv}(\xi) =$ $|\det(\dm)| \hat{\tv}(\dm^T \xi)
\hat{\ta}(\xi)$ and $ \widehat{\tz_{\ta,\dm}
\tv}(\dm^T \xi) =|\det(\dm)|^{-1} \sum_{\omega \in \dmfcg{\dm}}
\hat{\tv}(\xi + 2 \pi \omega) \overline{\hat{\ta}(\xi + 2 \pi
\omega)}$, where $\hat{\ta}(\xi):=\sum_{\vk\in \ZZ^d} a(\vk) e^{-\iu \vk\cdot
\xi}$ and $\dmfcg{\dm}:=[(\dm^T)^{-1}\ZZ^d] \cap [0,1)^d$.

We will now illustrate the general idea of an AMRA to the reader in the
case of shearlet systems.


\subsection{Shearlet Systems}
\label{subsec:shearlets}

\begin{figure}[ht]
\begin{center}
\includegraphics[height=1.4in]{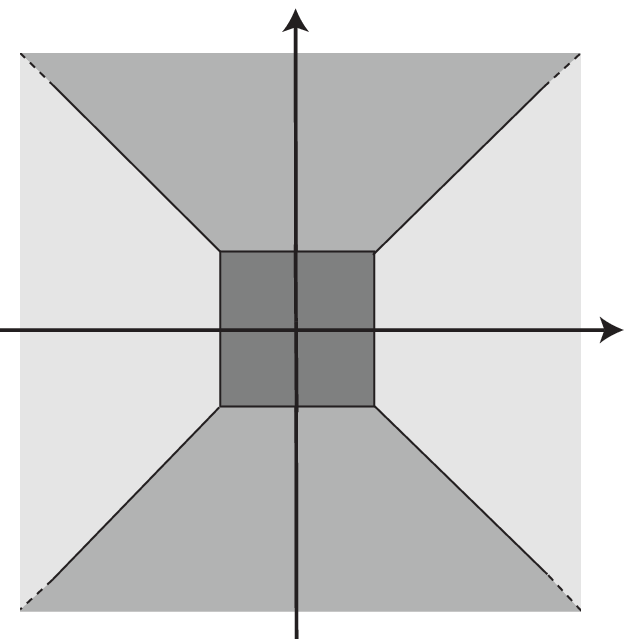}
\put(-33,58){\footnotesize{$\cC_0$}}
\put(-70,80){\footnotesize{$\cC_1$}}
\put(-88,30){\footnotesize{$\cC_0$}}
\put(-50,52){\footnotesize{$\cR$}}
\put(-45,15){\footnotesize{$\cC_1$}} \put(-60,-17){(a)}
\hspace*{4cm}
\includegraphics[height=1.4in]{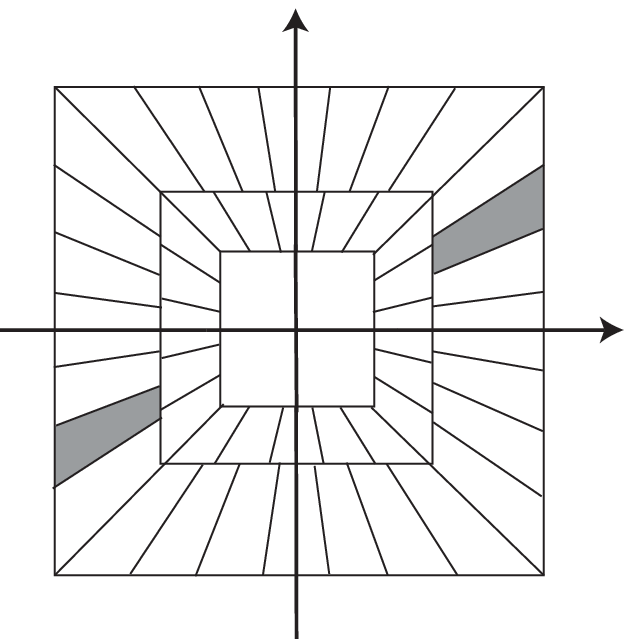}
\put(-60,-17){(b)}
\end{center}
\caption{(a) The cones $\cC_0$ and $\cC_1$ and the centered rectangle
$\cR$ in frequency domain. (b) The tiling of the frequency domain
induced by discrete shearlets.} \label{fig:ShearletsCone}
\end{figure}

Shearlet systems can be regarded as a special case of the previously
introduced affine-like systems. The continuum shearlet transform
with discrete parameters \cite{GKL06,KL09} for functions in $L^2(\bR^2)$
uses a two-parameter dilation group, where one parameter indexes scale,
and the second parameter indexes orientation. For each $c>0$ and $s \in
\mathbb{R}$, let $A_c$ denote the {\em parabolic scaling matrix} and
$S_s$ denote the {\em shear matrix} of the form
\[
A_c =
\left(\begin{array}{cc}
  c & 0\\ 0 & \sqrt{c}
\end{array}\right)
\qquad\mbox{and}\qquad S_s = \left(\begin{array}{cc}
  1 & s\\ 0 & 1
\end{array}\right),
\]
respectively. To provide an equal treatment of the $x$- and
$y$-axis, we split the frequency plane into the horizontal cone
\[
\cC_0 = \{(\xi_1,\xi_2) \in \bR^2 : |\xi_1| \ge 1,\, |\xi_1/\xi_2| \ge 1\},
\]
the vertical cone
\[
\cC_1 = \{(\xi_1,\xi_2) \in \bR^2 : |\xi_2| \ge 1,\, |\xi_1/\xi_2| \le 1\},
\]
as well as a centered rectangle
\[
\cR = \{(\xi_1,\xi_2) \in \bR^2 : \|(\xi_1,\xi_2)\|_\infty < 1\}
\]
(see Figure \ref{fig:ShearletsCone} (a)).

For cone $\cC_1$, at scale $j \ge 0$, orientation $k = -2^{j}, \dots,$
$2^{j}$, and spatial position $m \in \bZ^2$,  the associated {\em shearlets}
are then defined by their Fourier transforms
\[
\sigma_{\eta} = 2^{j\frac{3}{4}} \psi(S_k A_{4^{j}} \cdot\,-m),
\]
where $\eta = (j,k,m,\iota)$ index scale, orientation, position, and
cone. The shearlets for $\cC_1$ are defined likewise by symmetry, as
illustrated in Figure \ref{fig:ShearletsCone} (b), and we denote the
resulting {\em shearlet system} by
\beq \label{eq:shearletsystem}
\{\sigma_{\eta} : \eta \in \bN_0 \times
\{-2^{j}, \dots, 2^{j}\} \times \bZ^2 \times \{0,1\} \}.
\eeq
This is an affine-like system as defined before.

Notice that we chose a scaling of $4^j$. Shearlet systems can be defined
similarly for a scaling of $2^j$, however, in this case the odd scales
have to be handled particularly carefully.

One particular interesting shearlet generator is the function $\psi \in L^2(\bR^2)$
defined by
\[
\hat{\psi}(\xi) = \hat{\psi}(\xi_1,\xi_2) = \hat{\psi}_1(\xi_1) \,
\hat{\psi}_2(\tfrac{\xi_2}{\xi_1}),
\]
where $\psi_1 \in L^2(\bR)$ is a wavelet with $\hat{\psi}_1 \in
C^\infty(\mathbb{R})$ and supp $\hat{\psi}_1 \subseteq [-4,-\frac14]
\cup [\frac14,4]$, and $\psi_2 \in L^2(\mathbb{R})$ is a `bump'
function satisfying $\hat{\psi}_2 \in C^\infty(\mathbb{R})$ and supp
$\hat{\psi}_2 \subseteq [-1,1]$. Filling in the low frequency band
appropriately, with this particular generator the shearlet system
\eqref{eq:shearletsystem} can be proven (see \cite[Thm. 3]{GKL06})
to form a tight frame for $\{f \in L^2(\bR^2) : \Sp \hat{f} \subset \cC_0 \cup \cC_1\}$.
We remind the reader that a system $X$ is
called a {\em tight frame} for a Hilbert space $\mathcal{H}$
(in the literature sometimes also referred to as a {\em Parseval
frame} or a {\em tight frame with bound one}), if
$\|f\|^2=\sum_{g\in X}|\langle f, g\rangle|^2$ holds for all $f\in
\mathcal{H}$.

Concluding, the definition just discussed shows that shearlets live on anisotropic regions of
width $2^{-2j}$ and length $2^{-j}$ at various orientations, which
are parametrized by slope rather than angle as for second generation
curvelets.


\subsection{Shearlet AMRA}
\label{subsec:shearletAMRA}

An intriguing application of the theory whose core ideas were
presented in the Subsection \ref{subsec:AMRA} is the consideration of
shearlet systems. In fact, this allows us to derive a framework
which fulfills the desiderata we aimed for (see Subsection
\ref{subsec:desiderata}).

The fast spatial decomposition algorithm is derived as a special
case of the methodology we will develop for general affine-like systems.
To illustrate the idea of AMRA in the situation of shearlet
systems, we will discuss one step of the decomposition algorithm.
For this, let $S_\mi$, $\mi = 1,\ldots,r$ be a selection of $2
\times 2$ shear matrices, i.e., matrices of the type $S_k$, $k \in
\ZZ$, let $1 \le s \le r$ be the separator between low- and
high-frequency part, and let $a_\mi$, $\mi = 1, \ldots, \mpsi$ be
finitely supported masks. Further, we denote the given data by $\tv$
which for convenience purposes we now assume to lie in
$\llll(\ZZ^2)$. Notice that certainly from the previous
decomposition step we presumably have many such low-frequency
coefficients. Retaining the notion introduced in Subsection \ref{subsec:AMRA},
we then compute the next level of low-frequency
coefficients by
\[
\tv_\ell = \tz_{\ta_\ell, S_{\ell}A_4} \tv,\qquad \ell = 1, \ldots,
s.
\]
Similarly, we compute the next level of high-frequency coefficients
by
\[
\tv_\ell = \tz_{\ta_\ell, S_{\ell}A_4} \tv,\qquad \ell = s+1,
\ldots, r.
\]
The next step continues with decomposing $\tv_\ell$,
$\ell=1,\ldots,s$. The total number of decomposition steps is $J$,
hence finite.

The reader should notice the requirements (R1) and (R2) from
Subsection \ref{subsec:AMRA} as opposed to a `classical MRA-decomposition
algorithm'. (R1) is the finite number of decomposition steps, which
is a more accurate model for what truly happens in a numerical
decomposition. It thereby allows more flexibility, since the limit
$J \to \infty$ doesn't need to be taken into account.
(R2) can be seen by the fact that the separator between low- and
high-frequency part is somehow `loose' in the sense that at each
level a certain condition needs to be satisfied (see  Theorem
\ref{theo:generalshearletUEP} (i) below), which does {\em not} distinguish
between those parts. This also implies that the structure of the
subspaces the data is projected onto is not that strict than for a
classical MRA, but allows also non-orthogonality and
non-inclusiveness.

This algorithm is accompanied by a perfect reconstruction algorithm,
if and only if, in each decomposition step the following version of
the Unitary Extension Principle \cite{RS97a} is satisfied by the
masks $\ta_{\mi}$. For this, we have the following result, which we
could coin the `Shearlet Unitary Extension Principle'.

\begin{theorem}
\label{theo:generalshearletUEP} Let $S_\mi$, $\mi = 1,\ldots,r$ be a
selection of $2 \times 2$ shear matrices, i.e., matrices of the type
$S_k$, $k \in \ZZ$, and let $\mi = 1, \ldots, \mpsi$, be finitely
supported masks. Then the following conditions are equivalent.
\bitem
\item[{\rm (i)}] For all $\tv \in \llll(\ZZ^d)$,
\[
\sum_{\mi=1}^\mpsi \sd_{\ta_\mi,S_\mi A_4} \tz_{\ta_\mi,S_\mi A_4}
\tv = \tv.
\]
\item[{\rm (ii)}] For any $\omega \in \Omega = \bigcup_{\mi=1}^\mpsi \dmfcg{S_\mi A_4}$, where $\dmfcg{S_\mi A_4}
:=[(S_\mi^T)^{-1} (\tfrac14\ZZ \times \tfrac12\ZZ)]\cap [0,1)^2$,
\[
\sum_{\mi \in \{n\; :\; \omega \in \dmfcg{S_n A_4}\}}
\widehat{\ta_\mi}(\xi) \overline{\widehat{\ta_\mi}(\xi + 2 \pi
\omega)} = \delta(\omega),
\]
where $\delta$ denotes the Dirac sequence such that $\delta(0)=1$
and $\delta(\omega)=0$ for $\omega \ne 0$. \eitem

\end{theorem}

In Section \ref{sec:general}, we state a general version of this
result (Theorem \ref{theo:generalUEP}) from which this theorem can
be derived as a corollary. A continuum version of it will then be
discussed in Section \ref{sec:continuum}. In Section
\ref{sec:shearlets}, we will provide a variety of different choices
for the masks, the dilation matrices, and the direction matrices.
Adapted to the shearlet setting, the dilation matrices will
typically be parabolic scaling matrices and the direction matrices
will be chosen to be shear matrices.


\subsection{Extensions and Further Questions}

In fact, the results -- including the general affine-like systems we
claim to consider --  are susceptible of extensive generalizations
and extensions most of which are far beyond the scope of this paper
and will be studied in future work.

\bitem
\item {\em Bi-Frame Case.} The bi-orthogonal framework for wavelets typically
allows for significantly weakened filter conditions. Also for our
more general MRA-framework, we can ask a similar question, which is
then correctly stated as the bi-frame case.
\item {\em Shearlets for Higher Dimensions.}
Thriving applications such as the problem of geometric separation in
image processing, for example, in astronomy when images of
galaxies require separated analyses of stars, filaments, and sheets
call for directional representation systems for 3D data, but also
even higher dimensions (see, e.g.,  \cite{DK08a,DK08aa}). Our framework
conveys the potential to generate shearlets for higher dimensional
signals alongside a fast flexible decomposition strategy, and we
will present one example in 3D in Section \ref{sec:shearlets}.
\item {\em Computational Realization.} Certainly, one extension of the work presented
in this paper is the algorithmic study of the algorithms presented
in this paper. Here we aim for understanding the additional
flexibility these provide and how to optimally explore and utilize
this property. \eitem


\subsection{Related Work}

Several research teams have previously designed MRA decomposition
algorithms based on shearlets: we mention the affine system-based
approach \cite{GLLWW06}, the subdivision-based approach \cite{KS08},
and the approach based on separability \cite{Lim09}. However,
neither of these approaches did satisfy all items of our list of
desiderata (see Subsection \ref{subsec:desiderata}). Further non-MRA
based approaches were undertaken, for instance, in \cite{ELL08}. In
our opinion, these pioneer efforts demonstrate real progress in
directional representation, but further progress is needed to derive
an in-all-aspects satisfactory comprehensive study of a fast spatial
domain shearlet transform within an appropriate MRA framework with
careful attention to mathematical exactness, faithfulness to the
continuum transform, and computational feasibility, ideally
fulfilling all our desiderata.

A particular credit deserves the work in \cite{KS08}, in which the
adaptivity ideas were already lurking. The main difference to this
paper is the additional freedom provided by the AMRA structure.


\subsection{Contribution of this Paper}

The contribution of this paper is three-fold. First, we introduce
the concept of an adaptive MRA (AMRA) structure suited to the main goal of
a fast flexible decomposition strategy. Secondly, we study this
novel methodology for the general case of affine-like systems. And
thirdly, we present a comprehensive theory for shearlet systems
which encompasses tight shearlet frames with spatially compactly
supported generators within such an AMRA structure along
with a fast decomposition strategy.


\subsection{Contents}

In Section \ref{sec:general}, we introduce the notation we employ
for general affine-like systems, state the fast decomposition algorithm
based on an AMRA, and prove the Unitary Extension Principle
(UEP) for this situation. Section \ref{sec:continuum} is
concerned with the relation to the continuum setting, i.e.,
with developing characterizing equations and approximation
properties for the functions associated with an ARMA. This general
methodology is then applied to the situation of shearlet systems in
Section \ref{sec:shearlets}, where we present a general construction
for associated filters. Here we also consider shearlet-like systems
in the sense of integer-valued matrices approximating rotations
different from the customarily employed shear matrices, and an
approach of detecting anisotropic phenomena with an astonishingly simple
isotropic system.


\section{An Adaptive Multiresolution Analysis for General Affine-Like Systems}
\label{sec:general}

As already elaborated upon in the introduction, one main idea of an
AMRA is to be able to design each decomposition step adaptively, for
instance, dependent on the previous decomposition. To follow this
philosophy, in Subsection \ref{subsec:UEP}, we will firstly analyze
one single decomposition step, which might occur at any stage of the
general decomposition algorithm. Secondly, in Subsection
\ref{subsec:fastdecomp}, we will then present the large picture in
the sense of the complete decomposition procedure.

\subsection{A Unitary Extension Principle for One Decomposition Step}
\label{subsec:UEP}

Let now $\tv \in \llll(\ZZ^d)$ be some set of data. This could be
the initial data, but also data after some steps of decomposition
then on a renormalized grid. We assume that we are given a sequence
of arbitrary $d \times d$ matrices $\dm_\mi$, $1 \le \mi \le r$, and
finitely supported masks $\ta_\mi$, $\mi = 1, \ldots, r$ according
to which the data shall be decomposed. Our first result is a Unitary Extension Principle (UEP) for this situation, which characterizes
those matrices and filters, which allow perfect reconstruction from
the decomposed data using subdivision.

\begin{theorem}\label{theo:generalUEP}
Let $\dm_\mi$, $1 \le \mi \le \mpsi$ be $d \times d$ invertible
integer matrices, and let $\ta_\mi$, $\mi = 1, \ldots, \mpsi$, be
finitely supported masks. Then the following conditions are
equivalent. \bitem
\item[{\rm (i)}] For all $\tv \in \llll(\ZZ^d)$,
\[
\sum_{\mi=1}^\mpsi \sd_{\ta_\mi,\dm_\mi} \tz_{\ta_\mi,\dm_\mi} \tv =
\tv.
\]
\item[{\rm (ii)}] For any $\omega \in \Omega = \bigcup_{\mi=1}^\mpsi \dmfcg{\dm_\mi}$, where $\dmfcg{\dm_\mi}
:=[(\dm_\mi^T)^{-1}\ZZ^d]\cap [0,1)^d$,
\begin{equation}\label{filter:general}
\sum_{\mi \in \{n\; :\; \omega \in \dmfcg{\dm_n}\}}
\widehat{\ta_\mi}(\xi) \overline{\widehat{\ta_\mi}(\xi + 2 \pi
\omega)} = \delta(\omega).
\end{equation}
\eitem
\end{theorem}

\begin{proof}
First, notice that, for all $\tv \in \llll(\ZZ^d)$, (i) is
equivalent to
\begin{eqnarray*}
\hat{\tv}(\xi) & = & \sum_{\mi=1}^\mpsi \widehat{[\sd_{\ta_\mi,\dm_\mi} \tz_{\ta_\mi,\dm_\mi} \tv]}\\
& = & \sum_{\mi=1}^\mpsi |\det(\dm_\mi)| \widehat{\tz_{\ta_\mi,\dm_\mi} \tv} (\dm_\mi^T \xi) \widehat{\ta_\mi}(\xi)\\
& = & \sum_{\mi=1}^\mpsi \sum_{\omega \in \dmfcg{\dm_\mi}}
\hat{\tv}(\xi + 2 \pi \omega) \widehat{\ta_\mi}(\xi)
\overline{\widehat{\ta_\mi}(\xi + 2 \pi \omega)}.
\end{eqnarray*}
Employing the equivalence relation given by $\dmfcg{\dm_n}$, we can
rewrite the previous equation as \beq \label{eq:equivalent(i)}
\hat{\tv}(\xi) = \sum_{\omega \in \Omega} \sum_{\mi \in \{n\; :\;
\omega \in \dmfcg{\dm_n}\}} \hat{\tv}(\xi + 2 \pi \omega)
\widehat{\ta_\mi}(\xi) \overline{\widehat{\ta_\mi}(\xi + 2 \pi
\omega)} \quad \mbox{for all }\; \tv \in \llll(\ZZ^d). \eeq Since
(ii) implies \eqref{eq:equivalent(i)}, this proves (ii)
$\Rightarrow$ (i).

We now assume that (i), hence \eqref{eq:equivalent(i)} holds. We aim
to prove that this implies (ii). For this, we let
$B_\epsilon(\xi_0)$ denote as usual an open ball around $\xi_0$ with
radius $\epsilon$. We now first observe that for any $\omega_0 \in
\Omega$ and any $\xi_0 \in \RR$, there exists some $\tv \in
\llll(\ZZ^d)$ and $\epsilon > 0$ such that \bitem
\item[(a)] $\hat{\tv}(\xi + 2 \pi \omega_0) = 1$ for all $\xi \in B_\epsilon(\xi_0)$,
\item[(b)] $\hat{\tv}(\xi + 2 \pi \omega) = 0$ for all $\xi \in B_\epsilon(\xi_0), \omega \in \Omega \setminus \{\omega_0\}$,
\item[(c)] $\Sp\, \hat{\tv} \subseteq 2 \pi \omega_0 + B_{2 \epsilon}(\xi_0)$,
\eitem simply since $\Omega$ is discrete. Now we can conclude that
\eqref{eq:equivalent(i)} implies
\[
\hat{\tv}(\xi) = \sum_{\mi \in \{n \; :\; \omega \in
\dmfcg{\dm_n}\}} \widehat{\ta_\mi}(\xi)
\overline{\widehat{\ta_\mi}(\xi + 2 \pi \omega)} \quad \mbox{for all
}\, \omega \in \Omega,\,\xi \in B_\epsilon(\xi_0).
\]
Hence again by (a) -- (c), condition (ii) follows for all $\xi \in
B_\epsilon(\xi_0)$. Since $\xi_0$ is arbitrary chosen, condition
(ii) follows.
\end{proof}

The first part (equivalence of conditions (i) and (ii)) of the next
result follows immediately from Theorem \ref{theo:generalUEP} as a
special case. It significantly simplifies the condition on the mask
imposed by the UEP, provided the situation allows to always use the
same sampling lattice. Condition (iii) is the spatial domain
expression for the UEP condition (ii), which illustrates the filter
design problem in spatial domain. Certainly, this condition could
also be stated in the general situation of Theorem
\ref{theo:generalUEP}. We however decide to omit this, since the
content would be clouded by very technical details.

\begin{corollary} \label{coro:generalUEPspecial}
Let $\dm_\mi$, $1 \le \mi \le \mpsi$ be $d \times d$ invertible
integer matrices satisfying that $\dm_\mi \ZZ^d = \dm \ZZ^d$ for all
$1 \le \mi \le \mpsi$ for some matrix $\dm$, and let $\ta_\mi$, $\mi
= 1, \ldots, \mpsi$, be finitely supported masks. Then the following
conditions are equivalent. \bitem
\item[{\rm (i)}] For all $\tv \in \llll(\ZZ^d)$,
\[
\sum_{\mi=1}^\mpsi \sd_{\ta_\mi,\dm_\mi} \tz_{\ta_\mi,\dm_\mi} \tv =
\tv.
\]
\item[{\rm (ii)}] For any $\omega \in \dmfcg{\dm} = [(\dm^T)^{-1}\ZZ^d]\cap[0,1)^d$,
\[
\sum_{\mi =1}^\mpsi \widehat{\ta_\mi}(\xi)
\overline{\widehat{\ta_\mi}(\xi + 2 \pi \omega)} = \delta(\omega).
\]
\item[{\rm (iii)}] For all $\vk, \gamma \in \ZZ^d$,
\[
\sum_{\vn \in \ZZ^2} \overline{\ta_\mi(\vk+\dm\vn+\gamma)}
\ta_\mi(\dm \vn+\gamma) = |\det(\dm)|^{-1} \delta(\vk).
\]
\eitem
\end{corollary}

\begin{proof}
As already mentioned, the equivalence of (i) and (ii) does follow
from Theorem \ref{theo:generalUEP}.

We next prove equivalence between (ii) and (iii). By using the
definition of $\widehat{\ta_\mi}$, condition (ii) is equivalent to
\begin{eqnarray}\nonumber
\delta(\omega) & = & \sum_{\mi =1}^\mpsi \sum_{\vk \in \ZZ^d} \overline{\ta_\mi(\vk)}e^{\iu
\vk\cdot\xi}  \sum_{\vn \in \ZZ^d} \ta_\mi(\vn)e^{-\iu\vn\cdot(\xi +
2 \pi \omega)}\\ \label{eq:equivalence(iii)}
& = & \sum_{\mi =1}^\mpsi \sum_{\vk, \vn \in \ZZ^d}
\overline{\ta_\mi(\vk)}\ta_\mi(\vn) e^{\iu(\vk-\vn)\cdot\xi}
e^{-\iu \vn \cdot 2 \pi \omega}.
\end{eqnarray}
Next we denote $\dmcg{\dm} = \ZZ^d/[\dm\ZZ^d]$ with $0 \in \dmcg{\dm}$. Then
\[
\ZZ^d = \dmcg{\dm} + \dm \ZZ^d
\]
follows directly. This now allows us to rewrite
\eqref{eq:equivalence(iii)} as
\begin{eqnarray*}
\delta(\omega) & = & \sum_{\mi =1}^\mpsi \sum_{\gamma \in
\dmcg{\dm}} \sum_{\vk, \vn \in \ZZ^d}
\overline{\ta_\mi(\vk)}\ta_\mi(\dm \vn+\gamma)
e^{\iu (\vk-\dm \vn - \gamma)\cdot \xi}  e^{-\iu (\dm \vn + \gamma) 2 \pi \omega}\\
& = & \sum_{\mi =1}^\mpsi \sum_{\gamma \in \dmcg{\dm}} \sum_{\vk,
\vn \in \ZZ^d} \overline{\ta_\mi(\vk + \dm \vn + \gamma)}\ta_\mi(\dm
\vn+\gamma) e^{\iu \vk\cdot \xi}  e^{-\iu \gamma\cdot 2 \pi \omega}.
\end{eqnarray*}
Considering this equation in matrix form
\[
\left(e^{-\iu \gamma\cdot 2 \pi \omega}\right)_{\omega \in
\dmfcg{\dm}, \gamma\in \dmcg{\dm}} \left(\sum_{\vk, \vn \in \ZZ^d}
\overline{\ta_\mi(\vk + \dm \vn + \gamma)}\ta_\mi(\dm
\vn+\gamma)e^{\iu \vk\cdot \xi}\right)_{\gamma \in \dmcg{\dm}} =
\left(\delta(\omega)\right)_{\omega \in \dmfcg{\dm}},
\]
we can conclude by taking the inverse that
\begin{eqnarray*}
\lefteqn{\left(\sum_{\vk, \vn \in \ZZ^d} \overline{\ta_\mi(\vk + \dm \vn + \gamma)}\ta_\mi(\dm \vn+\gamma)e^{\iu \vk\cdot \xi}\right)_{\gamma \in \dmcg{\dm}}}\\
& = & |\det(\dm_1)|^{-1}\left(e^{\iu \gamma\cdot 2 \pi
\omega}\right)_{\omega \in \dmfcg{\dm}, \gamma \in \dmcg{\dm}}
\left(\delta(\omega)\right)_{\omega \in \dmfcg{\dm}}\\
& = & |\det(\dm)|^{-1} \left(\begin{array}{c}1\\ \vdots \\ 1
\end{array}\right).
\end{eqnarray*}
Thus (ii) is equivalent to the equation
\[
\sum_{\vk, \vn \in \ZZ^d} \overline{\ta_\mi(\vk + \dm \vn +
\gamma)}\ta_\mi(\dm \vn+\gamma)e^{\iu \vk\cdot\xi} =
|\det(\dm)|^{-1} \quad \mbox{for all }\; \gamma \in \ZZ^d,
\]
which in turn is equivalent to (iii).
\end{proof}

This corollary implies, by defining the equivalence relation
$\sim$ on the $d \times d$ invertible integer matrices by
\[
\dm \sim \dm'\quad :\Longleftrightarrow \quad \dm\ZZ^d = \dm'\ZZ^d,
\]
 that only one representative of each involved class needs to satisfy
\eqref{filter:general}. In other words, only the generated lattices
$\dm\ZZ^d$ do matter in the condition equivalent to perfect
reconstruction. Inside each class we have very much freedom to
choose the dilation matrices as necessary by the application.
We refer the reader to \cite{Han:laa} for the application of this observation on the construction of wavelets.

\subsection{Fast Decomposition Algorithm}
\label{subsec:fastdecomp}

In the previous subsection, we derived characterizing conditions for
sequences of arbitrary $d \times d$ matrices $\dm_\mi$, $1 \le \mi
\le \mpsi$, and finitely supported masks $\ta_\mi$, $\mi = 1,
\ldots, \mpsi$ which allow perfect reconstruction in each step of
an AMRA. This now enables us to
present a general adaptive decomposition algorithm, where the
matrices and masks in each step can be chosen according to Theorem
\ref{theo:generalUEP}.

\begin{figure}[ht]
\centering
\includegraphics[height=2.7in]{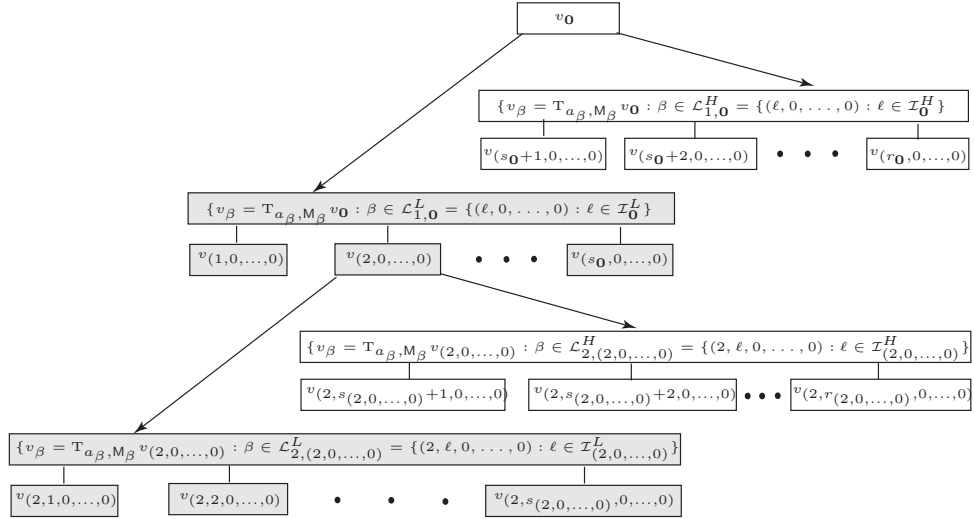}
\put(-158,188){{\tiny $v_\mathbf{0}$}}
\put(-180,154){{\tiny $\{v_\beta=\tz_{\ta_\beta, \dm_{\beta}} \tv_{\mathbf{0}} : \beta\in\cL^H_{1,\mathbf{0}} =\{(\ell,0,\ldots,0) : \ell \in \cI_{\mathbf{0}}^H\}$}}
\put(-185,137){{\tiny $v_{(s_\mathbf{0}+1,0,\ldots,0)}$}}
\put(-130,137){{\tiny $v_{(s_\mathbf{0}+2,0,\ldots,0)}$}}
\put(-39,137){{\tiny $v_{(r_\mathbf{0},0,\ldots,0)}$}}
\put(-290,115){{\tiny $\{v_\beta=\tz_{\ta_\beta, \dm_{\beta}} \tv_{\mathbf{0}} : \beta\in\cL^L_{1,\mathbf{0}} =\{(\ell,0,\ldots,0) : \ell \in \cI_{\mathbf{0}}^L\}$}}
\put(-293,98){{\tiny $v_{(1,0,\ldots,0)}$}}
\put(-237,98){{\tiny $v_{(2,0,\ldots,0)}$}}
\put(-152,98){{\tiny $v_{(s_\mathbf{0},0,\ldots,0)}$}}
\put(-252,63){{\tiny $\{v_\beta=\tz_{\ta_\beta, \dm_{\beta}} \tv_{(2,0,\ldots,0)} : \beta\in\cL^H_{2,(2,0,\ldots,0)} =\{(2,\ell,0,\ldots,0) : \ell \in \cI_{(2,0,\ldots,0)}^H\}$}}
\put(-252,47){{\tiny $v_{(2,s_{(2,0,\ldots,0)}+1,0,\ldots,0)}$}}
\put(-166,47){{\tiny $v_{(2,s_{(2,0,\ldots,0)}+2,0,\ldots,0)}$}}
\put(-68,47){{\tiny $v_{(2,r_{(2,0,\ldots,0)},0,\ldots,0)}$}}
\put(-362,25){{\tiny $\{v_\beta=\tz_{\ta_\beta, \dm_{\beta}} \tv_{(2,0,\ldots,0)} : \beta\in\cL^L_{2,(2,0,\ldots,0)} =\{(2,\ell,0,\ldots,0) : \ell \in \cI_{(2,0,\ldots,0)}^L\}$}}
\put(-362,7){{\tiny $v_{(2,1,0,\ldots,0)}$}}
\put(-300,7){{\tiny $v_{(2,2,0,\ldots,0)}$}}
\put(-181,7){{\tiny $v_{(2,s_{(2,0,\ldots,0)},0,\ldots,0)}$}}
\caption{The decomposition structure of {\sc (FAD)} illustrated
through the first complete step and the second step shown
exemplarily through the decomposition of $v_{(2,0,\ldots,0)}$. The low-frequency
components, which will be processed further, are gray-shaded. }
\label{fig:tree}
\end{figure}

We first require some notation to carefully keep track of the
decomposition steps and positions in the generated tree structure by
suitable indexes. While reading the definitions, we recommend the
reader to also take a look at Figure \ref{fig:tree}, which
illustrates the indexing of the decomposition. The heart of our
indices are vectors $(\beta_1, \ldots, \beta_J)$ assigned to
matrices, filters, and data, where each entry $\beta_j$ indicates
whether this object is related to the computation of scale $j$ or a
coarser scale, and if yes, whether the data at scale $j$ to reach
this object was generated by a low-pass or high-pass filter.

Let us now be more specific. The original data is assigned the index
$\mathbf{0}:= (0, \ldots,0) \in \mathbb{N}_0^J$, and we set $\cL_0^L
= \{\mathbf{0}\} = \{(0, \ldots,0)\}$. For the indexing of the
low-pass filters, we define $\cI_{\mathbf{0}}^L = \{1,\ldots,
\mphi_{\mathbf{0}} \}$ and, for the high-pass filters,
$\cI_{\mathbf{0}}^H = \{\mphi_{\mathbf{0}}+1, \ldots,
\mpsi_{\mathbf{0}}\}$, where $\mphi_{\mathbf{0}} \le
\mpsi_{\mathbf{0}}$ and $\mpsi_{\mathbf{0}}$ is a positive integer.
In that sense $\mphi_{\mathbf{0}}$ {\em partitions} the filters into
low- and high-frequency filters. If $\mphi_{\mathbf{0}} =
\mpsi_{\mathbf{0}}$, then $\cI_{\mathbf{0}}^H = \emptyset$. Thus,
matrices and filters in the first step of the decomposition are
labeled by $(\mi,0,\ldots,0) \in \mathbb{N}_0^J$, $\mi \in
\cI_{\mathbf{0}}^L$ for the low-frequency objects, and
$(\mi,0,\ldots,0) \in \mathbb{N}_0^J$, $\mi \in \cI_{\mathbf{0}}^H$
for the high-frequency objects. For our convenience, we further
introduce the notation
\[
\cL_{1}^L = \cL_{1,\mathbf{0}}^L = \{(\mi,0,\ldots,0) \in \RR^J :
\mi \in \cI_{\mathbf{0}}^L\}
\]
and
\[
\cL_{1}^H =
\cL_{1,\mathbf{0}}^H = \{(\mi,0,\ldots,0) \in \RR^J : \mi \in
\cI_{\mathbf{0}}^H\}
\]
to denote the indices for the low- and high-frequency objects in the
first decomposition layer. The associated filters $a_{\beta}, \beta
\in \cL_1^L \cup \cL_1^H$ for the first decomposition step are
required to satisfy condition (ii) in Theorem
\ref{theo:generalUEP}, i.e.,
\[
\sum_{\beta \in \{\gamma \in \cL_1^L \cup \cL_1^H\; :\; \omega \in
\dmfcg{\dm_{\gamma}}\}} \widehat{\ta_\beta}(\xi)
\overline{\widehat{\ta_\beta}(\xi + 2 \pi \omega)} = \delta(\omega),
\quad \omega \in \Omega_{\mathbf{0}}=\Omega_{(0,\ldots,0)},
\]
where $\dm_{\gamma}$ are $d\times d$ invertible integer matrices,
$\dmfcg{\dm_{\gamma}}:=[(\dm_{\gamma}^T)^{-1}\ZZ^d]\cap [0,1)^d$,
and $\Omega_{(0,\ldots,0)}:=\cup_{\gamma\in \cL^L_1 \cup \cL^H_1}
\dmfcg{\dm_\gamma}$.

The matrices and filters which are used in the $j$th step will then
be labeled as follows. We first assume that in the $(j-1)$th step
the sets $\cL^L_{j-1}$ and $\cL^H_{j-1}$ were already constructed.
Then, for some $\beta=(\beta_1, \ldots, \beta_{j-1},0,\ldots, 0)\in
\cL^L_{j-1}$, we define
\[
\cI^L_{\beta}=\{1,\ldots, \mphi_\beta\} \quad \mbox{and} \quad
\cI^H_\beta=\{\mphi_\beta+1, \ldots, \mpsi_\beta\}\quad \mbox{with
}1\le \mphi_\beta\le \mpsi_\beta,
\]
as single labels in the new decomposition layer for the matrices and
filters. To take the whole tree structure into account, we now
define the set of new low-pass related indices arising from
$\beta=(\beta_1, \ldots, \beta_{j-1},0, \ldots, 0)$ by
\[
\cL^L_{j,\beta}=\{ (\beta_1, \ldots, \beta_{j-1}, \mi, 0,\ldots, 0)
\in \mathbb{N}_0^J : \mi \in \cI^L_{\beta}\}
\]
and the same for the set of all high-pass related indices,
\[
\cL^H_{j,\beta}=\{ (\beta_1, \ldots, \beta_{j-1}, \mi, 0,\ldots, 0)
\in \mathbb{N}_0^J : \mi \in \cI^H_{\beta}\}.
\]
We further set
\[
\cL^L_{j}=\cup_{\beta\in \cL^L_{j-1}} \cL^L_{j,\beta} \quad
\mbox{and} \quad \cL^H_{j}=\cup_{\beta\in \cL^L_{j-1}}
\cL^H_{j,\beta},
\]
the complete set of low-pass and high-pass indices in step $j$,
respectively. Now let $\alpha \in \cL_{j-1}^L$. To ensure perfect
reconstruction, according to Theorem \ref{theo:generalUEP}, the
associated filters $\ta_{\beta}$, $\beta \in \cL_{j,\alpha}^L \cup
\cL_{j,\alpha}^H$  must satisfy
\[
\sum_{\beta \in \{\gamma \in \cL_{j,\alpha}^L \cup
\cL_{j,\alpha}^H\; : \; \omega\in \dmfcg{\dm_{\gamma}}\}}
\widehat{\ta_{\beta}}(\xi)\overline{\widehat{\ta_{\beta}}(\xi+2\pi
\omega)}=\delta(\omega), \qquad \omega\in \Omega_{\alpha},
\]
where $\Omega_{\alpha}:=\cup_{\gamma \in \cL_{j,\alpha}^L \cup
\cL_{j,\alpha}^H} \dmfcg{\dm_{\gamma}}$.

Next we describe the general multi-level  decomposition algorithm
explicitly using the indexing that we just introduced. For
illustrative purposes, before presenting the complete algorithm, we
display the first decomposition step as well as one part of the
second decomposition step in Figure \ref{fig:tree}.

\begin{figure}[ht]
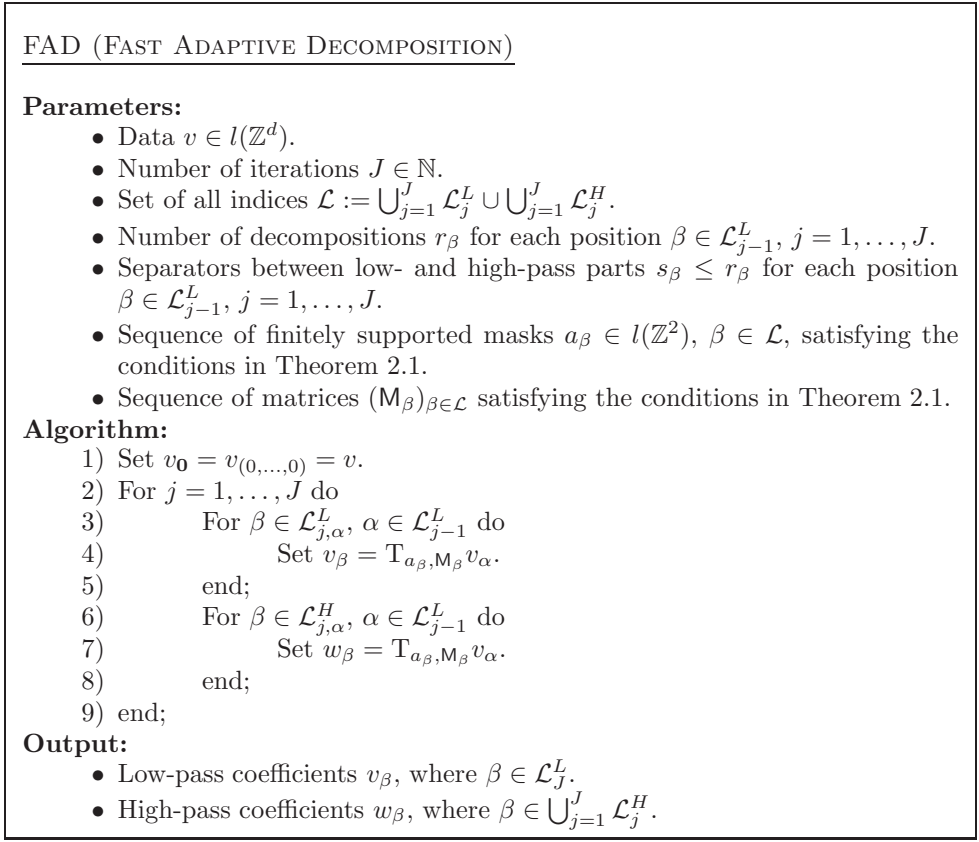

\centering \framebox{
\begin{minipage}[h]{4.9in}
\vspace*{0.3cm} {\sc \underline{FAD (Fast Adaptive Decomposition)}}

\vspace*{0.4cm}

{\bf Parameters:}\\[-3ex]
\begin{itemize}
\item Data $\tv \in \llll(\ZZ^d)$.
\item Number of iterations $J \in \NN$.
\item Set of all indices $\cL:=\bigcup_{j=1}^J \cL^L_{j} \cup \bigcup_{j=1}^J \cL^H_{j}$.
\item Number of decompositions $\mpsi_\beta$ for each position $\beta \in \cL^L_{j-1}$, $j=1,\ldots,J$.
\item Separators between low- and high-pass parts $\mphi_\beta \le \mpsi_\beta$
for each position $\beta \in \cL^L_{j-1}$, $j=1,\ldots,J$.
\item Sequence of finitely supported masks $\ta_{\beta} \in \llll(\ZZ^2)$, $\beta \in \cL$,
satisfying the conditions in Theorem \ref{theo:generalUEP}.
\item Sequence of matrices $(\dm_\beta)_{\beta \in \cL}$
satisfying the conditions in Theorem \ref{theo:generalUEP}.
\end{itemize}

{\bf Algorithm:}\\[-3ex]
\begin{itemize}
\item[1)] Set $v_{\mathbf{0}} = \tv_{(0,...,0)} = \tv$.
\item[2)] For $j=1, \ldots, J$ do
\item[3)] \hspace*{1cm} For $\beta \in \cL^L_{j,\alpha}$, $\alpha \in \cL^L_{j-1}$ do
\item[4)] \hspace*{2cm} Set $\tv_\beta = \tz_{\ta_\beta, \dm_{\beta}} \tv_{\alpha}$.
\item[5)] \hspace*{1cm} end;
\item[6)] \hspace*{1cm} For $\beta \in  \cL^H_{j,\alpha}$, $\alpha \in \cL^L_{j-1}$ do
\item[7)] \hspace*{2cm} Set $\tw_\beta = \tz_{\ta_\beta, \dm_{\beta}} \tv_{\alpha}$.
\item[8)] \hspace*{1cm} end;
\item[9)] end;
\end{itemize}

{\bf Output:}\\[-3ex]
\begin{itemize}
\item Low-pass coefficients
$\tv_{\beta}$, where $\beta \in \cL^L_{J}$.
\item High-pass coefficients
$\tw_{\beta}$, where $\beta \in \bigcup_{j=1}^J \cL^H_{j}$.
\end{itemize}
\vspace*{0.01cm}
\end{minipage}
} \caption{The {\sc FAD} Algorithm for a fast adaptive decomposition
using affine-like systems.} \label{fig:fad}
\end{figure}

In Figure \ref{fig:fad}, we now describe the general multi-level
decomposition algorithm explicitly. This decomposition can be
implemented as the usual fast wavelet transform with a tree
structure. We remind the reader of the introduced notion of
transition and subdivision in Subsection \ref{subsec:AMRA}.

\begin{figure}[th]
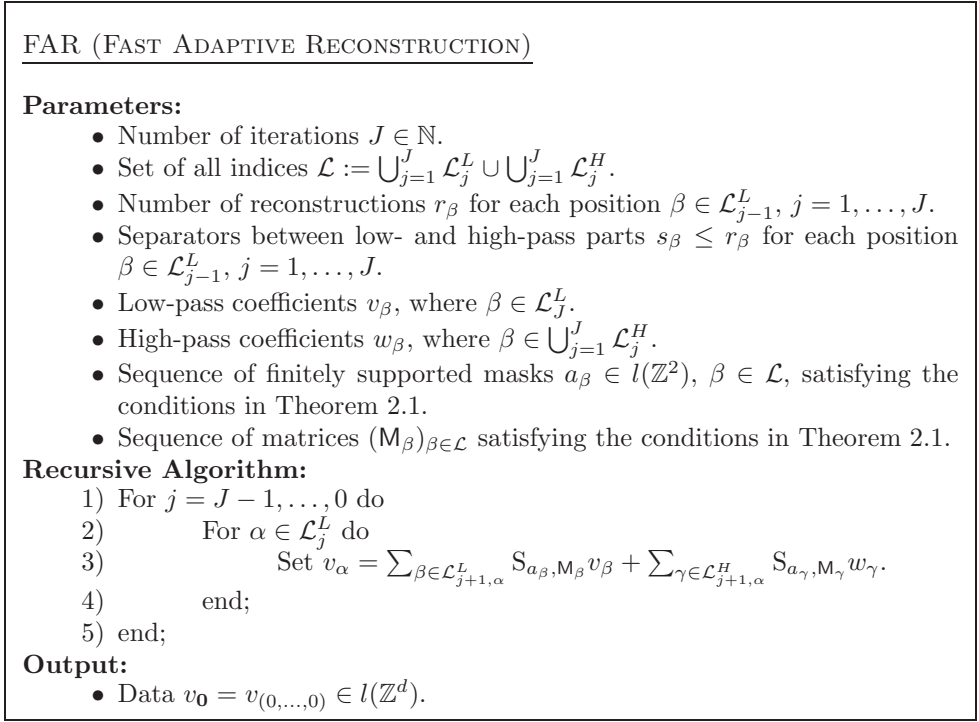

\centering \framebox{
\begin{minipage}[h]{4.9in}
\vspace*{0.3cm} {\sc \underline{FAR (Fast Adaptive Reconstruction)}}

\vspace*{0.4cm}

{\bf Parameters:}\\[-3ex]
\begin{itemize}
\item Number of iterations $J \in \NN$.
\item Set of all indices $\cL:=\bigcup_{j=1}^J \cL^L_{j} \cup \bigcup_{j=1}^J \cL^H_{j}$.
\item Number of reconstructions $\mpsi_\beta$ for each position $\beta \in \cL^L_{j-1}$, $j=1,\ldots,J$.
\item Separators between low- and high-pass parts $\mphi_\beta \le \mpsi_\beta$
for each position $\beta \in \cL^L_{j-1}$, $j=1,\ldots,J$.
\item Low-pass coefficients $\tv_{\beta}$, where $\beta \in \cL^L_{J}$.
\item High-pass coefficients $\tw_{\beta}$, where $\beta \in \bigcup_{j=1}^J \cL^H_{j}$.
\item Sequence of finitely supported masks $a_{\beta} \in \llll(\ZZ^2)$,  $\beta \in \cL$,
satisfying the conditions in Theorem \ref{theo:generalUEP}.
\item Sequence of matrices $(\dm_\beta)_{\beta \in \cL}$
satisfying the conditions in Theorem \ref{theo:generalUEP}.
\end{itemize}

{\bf Recursive Algorithm:}\\[-3ex]
\begin{itemize}
\item[1)] For $j=J-1, \ldots, 0$ do
\item[2)] \hspace*{1cm} For $\alpha \in \cL^L_{j}$ do
\item[3)] \hspace*{2cm} Set $\tv_{\alpha} = \sum_{\beta \in \cL_{j+1,\alpha}^L} \sd_{\ta_{\beta}, \dm_{\beta}} \tv_{\beta}+\sum_{\gamma \in \cL_{j+1,\alpha}^H} \sd_{\ta_{\gamma}, \dm_{\gamma}} \tw_{\gamma}$.
\item[4)] \hspace*{1cm} end;
\item[5)] end;
\end{itemize}

{\bf Output:}\\[-3ex]
\begin{itemize}
\item Data $v_{\mathbf{0}} = \tv_{(0,\ldots,0)} \in \llll(\ZZ^d)$.
\end{itemize}
\vspace*{0.01cm}
\end{minipage}
} \caption{The {\sc FAR} Algorithm for a fast adaptive
reconstruction using affine-like systems.} \label{fig:far}
\end{figure}

The fact that the filters are chosen to be perfect reconstruction
filters allows us to reconstruct the data by application of
appropriate subdivision operators as displayed in Figure
\ref{fig:far}.

\section{Relation with Continuum Setting}
\label{sec:continuum}

We next intend to relate the `digital conditions' on the filters to
`continuum conditions' for associated function systems, in
particular, frame systems.

\subsection{Characterization Equations}

As before in Theorem \ref{theo:generalUEP}, we only consider one
level. The conditions we will derive then need to be satisfied for
each step in the iteration, and the choice of a non-stationary or
stationary scheme is left to the user.

For a function $\phi$ and invertible $d\times d$ matrix $U$, we
shall use the notation as in \eqref{dilation:shift}.

We can now formulate the condition on prefect reconstruction in
terms of a condition in the function setting. Notice that here we
again consider the most general situation of Theorem
\ref{theo:generalUEP} as opposed to Corollary
\ref{coro:generalUEPspecial}.

\begin{theorem}\label{thm:onelevel}
Let $\dm_\mi$, $0 \le \mi \le \mpsi$ be $d \times d$ invertible
integer matrices, and let $\ta_\mi$, $\mi = 1, \ldots, \mpsi$, be
finitely supported masks. Let $\phi$ be a nontrivial compactly
supported function in $L^2(\RR^d)$. Then the following conditions
are equivalent. \bitem
\item[{\rm (i)}] For all $\tv \in \llll(\ZZ^d)$,
\[
\sum_{\mi=1}^\mpsi \sd_{\ta_\mi,\dm_\mi} \tz_{\ta_\mi,\dm_\mi} \tv =
\tv.
\]
\item[{\rm (ii)}] For each $f, g \in L^2(\RR^d)$,
\begin{equation}\label{identity:fg}
\sum_{\vk \in \ZZ^d} \langle f, \phi_{\dm_0;  \vk}\rangle \langle
\phi_{\dm_0; \vk}, g\rangle = \sum_{\mi = 1}^\mpsi \sum_{\vk \in
\ZZ^d} \langle f, \psi^\mi_{\dm_\mi^{-1} \dm_0;  \vk} \rangle
\langle \psi^\mi_{\dm_\mi^{-1}\dm_0;  \vk}, g\rangle,
\end{equation}
where $\phi_{\dm_0; \vk}$ is defined as in \eqref{dilation:shift}
and $\psi^1, \ldots, \psi^\mpsi$ are defined by
$\widehat{\psi^\mi}(\dm_\mi^T\xi):=\widehat{\ta_\mi}(\xi)\hat
\phi(\xi)$, that is,
\[
\psi^\mi := |\det (\dm_\mi)| \sum_{\vk\in\ZZ^d} \ta_\mi(\vk)
\phi(\dm_\mi \cdot - \vk).
\]
\eitem
\end{theorem}

\begin{proof}
We first recall that (i) is equivalent to \eqref{filter:general} by
Theorem \ref{theo:generalUEP}. The strategy will now be to show that
(ii) is equivalent to \eqref{filter:general}.

We first assume that \eqref{filter:general} holds. Then
$\langle f, \phi_{\dm_0; \vk}\rangle=\langle |\det \dm_0|^{-1/2}
f(\dm_0^{-1}\cdot), \phi(\cdot-\vk)\rangle=\langle |\det
\dm_0|^{-1/2} f(\dm_0^{-1}\cdot), \phi_{I_d; \vk}\rangle$. It is
easy to see that \eqref{identity:fg} holds if and only if it holds
for $\dm_0=I_d$. So, we assume that $\dm_0=I_d$ in the following
proof. Using the Fourier-based approach as in \cite[Lemma~3]{Han09},
we have
\[
\sum_{\vk \in \ZZ^d} \langle f, \phi_{I_d; \vk}\rangle \langle
\phi_{I_d; \vk}, g\rangle=(2\pi)^d \int_{\RR^d} \sum_{\vk\in \ZZ^d}
\hat f(\xi) \overline{\hat g(\xi+2\pi \vk)} \;
\overline{\hat{\phi}(\xi)} \hat \phi(\xi+2\pi \vk)d\xi
\]
and
\begin{equation}\label{func:eq2}
\begin{split}
&\sum_{\mi = 1}^\mpsi \sum_{\vk \in \ZZ^d} \langle f, \psi^\mi_{\dm_\mi^{-1}; \vk} \rangle \langle \psi^\mi_{\dm_\mi^{-1}; \vk}, g\rangle\\
&=(2\pi)^d \int_{\RR^d} \sum_{\mi=1}^\mpsi \sum_{\vk\in \ZZ^d} \hat
f(\xi)\overline{\hat g(\xi+2\pi (\dm_\mi^T)^{-1} \vk)} \;
\overline{\widehat{\psi^\mi}(\dm_\mi^T \xi)}
\widehat{\psi^\mi}(\dm_\mi^T \xi+2\pi \vk)d\xi.
\end{split}
\end{equation}
Since $\widehat{\psi^\mi}(\dm_\mi^T\xi):=\widehat{\ta_\mi}(\xi)\hat
\phi(\xi)$, we have
\[
\overline{\widehat{\psi^\mi}(\dm_\mi^T \xi)}
\widehat{\psi^\mi}(\dm_\mi^T \xi+2\pi \vk)
=\overline{\widehat{\ta_\mi}(\xi)} \widehat{\ta_\mi}(\xi+2\pi
(\dm_\mi^T)^{-1} \vk) \overline{\hat \phi(\xi)} \hat \phi(\xi+2\pi
(\dm_\mi^T)^{-1} \vk).
\]
Note that $\ZZ^d$ is the disjoint union of
$\dm_\mi^T\omega+\dm_\mi^T \ZZ^d$, $\omega\in \dmfcg{\dm_\mi}$. Now
we deduce that \eqref{func:eq2} becomes
\begin{align*}
&\sum_{\mi = 1}^\mpsi \sum_{\vk \in \ZZ^d} \langle f, \psi^\mi_{\dm_\mi^{-1} ; \vk} \rangle \langle \psi^\mi_{\dm_\mi^{-1}; \vk}, g\rangle\\
&=(2\pi)^d  \int_{\RR^d} \sum_{\mi=1}^\mpsi \sum_{\omega\in \dmfcg{\dm_\mi}} \sum_{\vk\in \ZZ^d} \hat f(\xi)\overline{\hat g(\xi+2\pi \omega+2\pi \vk)} \; \overline{\widehat{\ta_\mi}(\xi)} \widehat{\ta_\mi}(\xi+2\pi \omega)\\
&\cdot \overline{\hat \phi(\xi)}\hat\phi(\xi+2\pi \omega+2\pi \vk)d\xi\\
&=(2\pi)^d  \int_{\RR^d}  \sum_{\omega\in \Omega} \sum_{\vk\in
\ZZ^d} \hat f(\xi)\overline{\hat g(\xi+2\pi \omega+2\pi \vk)}
\;  \overline{\hat \phi(\xi)}\hat\phi(\xi+2\pi \omega+2\pi \vk)\\
& \cdot \sum_{\mi\in \{ 1\le n\le \mpsi\; : \;
\omega\in \dmfcg{\dm_n}\}} \overline{\widehat{\ta_\mi}(\xi)}
\widehat{\ta_\mi}(\xi+2\pi \omega)d\xi.
\end{align*}
Hence, \eqref{identity:fg} is equivalent to
\begin{eqnarray}\nonumber
\lefteqn{\int_{\RR^d} \sum_{\vk\in \ZZ^d} \hat f(\xi) \overline{\hat g(\xi+2\pi \vk)} \; \overline{\hat{\phi}(\xi)} \hat \phi(\xi+2\pi \vk)d\xi}\\ \nonumber
& = & \int_{\RR^d}  \sum_{\omega\in \Omega} \sum_{\vk\in \ZZ^d} \hat
f(\xi)\overline{\hat g(\xi+2\pi \omega+2\pi \vk)}
\;  \overline{\hat \phi(\xi)}\hat\phi(\xi+2\pi \omega+2\pi \vk)\\ \label{func:eq3}
& & \cdot \sum_{\mi \in \{ 1\le n\le \mpsi\; : \;
\omega\in \dmfcg{\dm_n}\}} \overline{\widehat{\ta_\mi}(\xi)}
\widehat{\ta_\mi}(\xi+2\pi \omega)d\xi.
\end{eqnarray}
Since \eqref{filter:general} holds, it is obvious that the above
identity in \eqref{func:eq3} holds and therefore, (iii) holds.

Conversely, if (ii) holds, then \eqref{func:eq3} holds. Now we use a
similar argument as in \cite[Lemma 5]{Han09} to prove that this implies condition (ii). For this, we
note that $\Omega+\ZZ^d:=\{ \omega+\vk\; : \; \omega\in \Omega,
\vk\in \ZZ^d\}$ is a discrete set without any accumulation point.
For any $\omega_0\in \Omega$ and $\xi_0\in \RR^d$, we can take any
functions $f,g \in L^2(\RR^d)$ such that the support of $\hat f$ is
contained inside $B_\epsilon(\xi_0)$ and the support of $\hat g$ is
contained inside $B_\epsilon(\xi+2\pi \omega_0)$. As long as
$\epsilon$ is small enough, it is not difficult to see that
\[
\hat f(\xi) \overline{\hat g(\xi+2\pi \vk)}=0 \qquad \mbox{for all }
\vk\in [\Omega+\ZZ^d] \backslash \{\omega_0\}.
\]
Now it is easy to see that \eqref{func:eq3} becomes
\begin{eqnarray*}
\lefteqn{\delta({\omega_0}) \int_{\RR^d} \hat f(\xi)\overline{\hat
g(\xi)}
\overline{\hat \phi(\xi)}\hat \phi(\xi)d\xi}\\
& = & \int_{\RR^d} \hat f(\xi)\overline{\hat g(\xi+2\pi \omega_0)}
\overline{\hat \phi(\xi)}\hat \phi(\xi+2\pi \omega_0)\cdot
\sum_{\mi\in \{ 1\le n\le \mpsi\; : \; \omega\in \dmfcg{\dm_n}\}}
\overline{\widehat{\ta_\mi}(\xi)} \widehat{\ta_\mi}(\xi+2\pi
\omega_0)d\xi.
\end{eqnarray*}
Note that $\hat \phi(\xi)\ne 0$ for almost every $\xi\in \RR$. Since
the above identity holds for all functions $f, g$ in $L^2(\RR^d)$ as
long as the support of $\hat f$ is contained inside
$B_\epsilon(\xi_0)$ and the support of $\hat g$ is contained inside
$B_\epsilon(\xi+2\pi \omega_0)$, we now easily deduce from the above
identity that $\sum_{\mi\in \{ 1\le n\le \mpsi\; : \; \omega\in
\dmfcg{\dm_n}\}} \overline{\widehat{\ta_\mi}(\xi)}
\widehat{\ta_\mi}(\xi+2\pi \omega_0) =\delta(\omega_0)$ for almost
every $\xi \in B_\epsilon(\xi_0)$. Since $\xi_0$ can be arbitrary,
we now conclude that \eqref{filter:general} holds.

The theorem is proved.
\end{proof}

 If $\mpsi=1$ and $\dm_1=I_d$ in Theorem \ref{thm:onelevel} (ii), then we must have $\widehat{\ta_1}=1$ (that is, $\ta_1=\delta$) and consequently for
this particular case, $\tz_{\ta_1, I_d} \tv=\tv$ and $\sd_{\ta_1,
I_d} \tv=\tv$, that is, the data is just copied.

\subsection{Tight Frame Structure}

For any positive integer $J$, we now construct a tight affine-like frame
$\cAS_J$ in $L^2(\RR^d)$ corresponding to the decomposition and
reconstruction algorithms, {\sc(FAD)} (Figure \ref{fig:fad}) and
{\sc(FAR)} (Figure \ref{fig:far}), respectively. Let $a$ be a mask
on $\ZZ^d$ and $M_0$ be a $d\times d$ dilation matrix, for example,
$\dm_0=2I_d$ and $a$ is a tensor product mask. We assume that
\[
\hat \phi(\xi):=\prod_{j=1}^\infty \hat\ta((\dm_0^T)^{-j}\xi),
\qquad \xi\in \RR^d
\]
is a well-defined function in $L^2(\RR^d)$ and we also assume that
there exist filters $\tb_1, \ldots, \tb_\mpsi$ such that
\[
\hat \ta(\xi)\overline{\hat \ta(\xi+2\pi \omega)}+\sum_{\mi=1}^\mpsi
\widehat{\tb_\mi}(\xi)\overline{\widehat{\tb_\mi}(\xi+2\pi \omega)}
=\delta(\omega), \qquad \omega\in \dmfcg{\dm_0}.
\]
Based on these, we then define
\[
\widehat{\psi^\mi}(\dm_0^T\xi):=\widehat{\tb_\mi}(\xi)\hat
\phi(\xi), \qquad \mi=1, \ldots, \mpsi.
\]
Then the system
$$
\mathcal{W}\mathcal{S}(\psi^1,\ldots, \psi^\mpsi):=\{ \psi^1_{\dm_0^j; \vk},
\ldots, \psi^\mpsi_{\dm_0^j; \vk}\; : \; j\in \ZZ, \vk\in \ZZ^d\}
$$
forms a tight frame in $L_2(\RR^d)$ by the Unitary Extension Principle of \cite{RS97a}
(see also \cite{HanShen:sima:08,HanShen:ca:09,RS97b}).

Next, we use this standard tight frame system to derive a new tight
affine-like frame $\cAS_J$ in $L^2(\RR^d)$ corresponding to the
decomposition and reconstruction algorithms, {\sc(FAD)} (Figure
\ref{fig:fad}) and {\sc(FAR)} (Figure \ref{fig:far}), respectively.

First, denote
\[
\psi^{(0,\ldots,0)}=\phi.
\]
Then, let $\beta \in \cL^L_{J} \cup \bigcup_{j=1}^J \cL^H_{j}$ with $\beta
\in \cL^L_{J,\beta^{(1)}}$, $\beta^{(1)} \in \cL^L_{J-1}$, or $\beta
\in \cL^H_{j,\beta^{(1)}}$, $\beta^{(1)} \in \cL^L_{j-1}$ for some
$j \in \{1, \ldots, J\}$. Then we recursively define affine
functions by
\[
\widehat{\psi^{\beta}}(\dm_{\beta}^T\xi):= \widehat{a_{\beta}}(\xi)
\widehat{\psi^{\beta^{(1)}}}(\xi).
\]
Further, let $\beta^{(2)} \in \cL^L_{j-2,\mu^{(2)}}$, etc. Using
this sequence, we recursively define matrices by
\[
\dn_\beta:=\dm_{\beta}^{-1} \dm_{\beta^{(1)}}^{-1} \cdots
\dm_{\beta^{(j-1)}}^{-1} \dm_{\beta^{(j)}}^{-1}.
\]

We now define the associated affine-like system in the following way:

\begin{definition} \label{defi:ASJ}
Retaining the introduced notions and definitions, the
{\em affine-like system $\cAS_J$} is defined by
\[
\cAS_J:=\{ \psi^1_{\dm_0^j; \vk}, \ldots, \psi^\mpsi_{\dm_0^j;
\vk}\; : \; j\ge J, \vk\in \ZZ^d\} \cup \{ \psi^\beta_{\dn_\beta
\dm_0^J; \vk}\; : \; \vk\in \ZZ^d, \beta \in \cL^L_{J} \cup
\bigcup_{j=1}^J \cL^H_{j}\},
\]
where we also employ the notation introduced in \eqref{dilation:shift}.
\end{definition}

We next show that this system -- although in general not forming an
orthonormal basis -- still always constitutes a tight frame.

\begin{theorem}\label{thm:twf}
For any positive integer $J$, $\cAS_J$ is a tight frame for
$L^2(\RR^d)$.
\end{theorem}

\begin{proof} By Theorem~\ref{thm:onelevel}, it is not difficult to deduce that
\[
\sum_{\beta\in \cL^L_{j} \cup \cL^H_{j}}\sum_{\vk\in \ZZ^d} |\langle
f, \psi^\beta_{\dn_\beta \dm_0^J; \vk}\rangle|^2 =\sum_{\beta\in
\cL^L_{j-1}}\sum_{\vk\in \ZZ^d} |\langle f, \psi^\beta_{\dn_\beta
\dm_0^J; \vk}\rangle|^2.
\]
Now from the above relation, we see that
\[
\sum_{\beta\in \cL^L_{J} \cup \bigcup_{j=1}^J \cL^H_{j}}\sum_{\vk\in
\ZZ^d} |\langle f, \psi^\beta_{\dn_\beta \dm_0^J; \vk}\rangle|^2
=\sum_{\vk\in \ZZ^d} |\langle f, \psi^{(0,\ldots,0)}_{\dm_0^J;
\vk}\rangle|^2 =\sum_{\vk\in \ZZ^d} |\langle f, \phi_{\dm_0^J;
k}\rangle|^2.
\]
By our given assumption on $\phi$ and $\psi^1, \ldots, \psi^\mpsi$,
it is known that $\{\phi_{\dm_0^J; \vk}\; : \; \vk\in \ZZ^d\}\cup
\{\psi^\mi_{\dm_0^j; \vk}\; : \; j\ge J, \vk\in \ZZ^d, \mi=1,
\ldots, \mpsi\}$ is a tight frame for $L^2(\RR^d)$. Consequently,
\[
\sum_{\beta \in \cL^L_{J} \cup \bigcup_{j=1}^J
\cL^H_{j}}\sum_{\vk\in \ZZ^d} |\langle f, \psi^\beta_{\dn_\beta
\dm_0^J; \vk}\rangle|^2 +\sum_{j=J}^\infty \sum_{\mi=1}^\mpsi
\sum_{\vk\in \ZZ^d} |\langle f, \psi^\mi_{\dm_0^j;
\vk}\rangle|^2=\|f\|^2_{L^2(\RR^d)}.
\]
Hence, $\cAS_J$ is a tight frame for $L^2(\RR^d)$.
\end{proof}

Next we analyze the approximation order of the affine-like system
$\cAS_J$. For this, for $\tau\ge 0$, we recall that $\dHH{\tau}$
consists of all functions $f \in L^2(\RR^d)$ satisfying
\[
|f|_{\dHH{\tau}}^2:=\frac{1}{(2\pi)^d} \int_{\RR^d} |\hat f(\xi)|^2
\|\xi\|^{2\nu} d\xi<\infty.
\]

We say that a mask $\ta: \dZ \mapsto \C$ has $\tau$ sum rules if
\[
\sum_{\vk\in \ZZ^d} \ta(\vn+\dm_0\vk) (\vn+\dm_0 \vk)^\beta
=\sum_{\vk\in \ZZ^d} \ta(\dm_0\vk) (\dm_0 \vk)^\beta
\]
for all $\vn\in \ZZ^d$ and for all $\beta=(\beta_1, \ldots,
\beta_d)\in (\NN \cup\{0\})^d$ such that $0\le \beta_1, \ldots,
\beta_d<\tau, \beta_1+\cdots+\beta_d<\tau$.

It follows from Theorem~\ref{thm:twf} that for any $f\in
L^2(\RR^d)$, expanding $f$ under the tight frame $\cAS_J$, we have
\[
f=\sum_{\beta \in \cL^L_{J} \cup \bigcup_{j=1}^J
\cL^H_{j}}\sum_{\vk\in \ZZ^d} \langle f, \psi^\beta_{\dn_\beta
\dm_0^J; \vk}\rangle \psi^\beta_{\dn_\beta \dm_0^J; \vk}
+\sum_{j=J}^\infty\sum_{\mi=1}^\mpsi \sum_{\vk\in \ZZ^d} \langle f,
\psi^\mi_{\dm_0^j; \vk}\rangle \psi^\mi_{\dm_0^j; \vk}
\]
with the series converging absolutely in $L^2(\RR^d)$. We now
consider the truncated series
\[
P_J f:=\sum_{\beta \in \cL^L_{J} \cup \bigcup_{j=1}^J
\cL^H_{j}}\sum_{\vk\in \ZZ^d} \langle f, \psi^\beta_{\dn_\beta
\dm_0^J; \vk}\rangle \psi^\beta_{\dn_\beta \dm_0^J; \vk}.
\]
By the same argument as in Theorem~\ref{thm:twf}, it is not
difficult to deduce that
\[
P_J f=\sum_{\vk\in \ZZ^d} \langle f, \phi_{\dm_0^J; \vk}\rangle
\phi_{\dm_0^J; \vk}
\]
The next theorem now follows from well-known results (see, e.g., \cite{DHRS,
Han:jat:01,Han:jat:03,Han:ACHA:09,Jia} and references therein).

\begin{theorem}\label{thm:approx}
Let $\cAS_J$ be the tight frame for $L^2(\RR^d)$ introduced in Definition \ref{defi:ASJ}.
Suppose that $\dm_0$ is isotropic, that is,
$\dm_0$ is similar to a diagonal matrix with all entries having the
same modulus. If the mask $\ta$ has $\tau$ sum rules, then $\cAS_J$
has  approximation order $\tau$, that is, there exists a positive
constant $C$ such that
\[
\|f-P_J f\|_{L^2(\RR^d)} \le C |\det \dm_0|^{-\tau J/d}
|f|_{\dHH{\tau}} \qquad \mbox{for all } f\in \dHH{\tau}, J\in \NN.
\]
\end{theorem}


\section{Shearlet Systems With an Associated AMRA}
\label{sec:shearlets}

In this section, we will apply our general framework to the special
case of shearlets.   This leads to shearlet systems, especially compactly
supported ones, associated with an AMRA structure and fast decomposition and
reconstruction algorithms. We anticipate that our considerations will
improve the applicability of shearlets to various applications, in particular,
frame based image restorations.

We first present a very general construction, which -- as shown in the
second subsection --  can be utilized to construct various shearlet
systems with an associated AMRA and fast decomposition and reconstruction
algorithms. Subsection \ref{subsec:3D} is then concerned with 3D shearlet
constructions.

\subsection{General Construction of an AMRA}
\label{subsec:generalconstr}

We first follow up on the comments after Corollary
\ref{coro:generalUEPspecial}, in which it was pointed out that the
construction of a sequence of dilation matrices $\dm_\mi$, $1 \le
\mi \le \mpsi$, i.e., invertible integer matrices and filters
$\ta_\mi$, $\mi = 1, \ldots, \mpsi$ satisfying
\eqref{filter:general}, depends only on the generated lattice
$\dm_\mi \ZZ^d$, $1 \le \mi \le \mpsi$. Hence we first present a
very general construction of such a sequence of matrices and filters
provided that the matrices satisfy $\dm_\mi \ZZ^d = \dm \ZZ^d$ for
all $1 \le \mi \le \mpsi$ for some matrix $\dm$, which was the
hypothesis of Corollary \ref{coro:generalUEPspecial}.

For this, we first fix a sublattice $\Gamma$ of $\ZZ^d$. Obviously,
there exist many $d\times d$ matrices $\dm$ such that $\dm \dZ
=\Gamma$. We now construct a set of tight frame filters for such a
lattice $\Gamma$ so that
\begin{equation}\label{cond:basic}
\sum_{\ell=1}^\mpsi
\widehat{\ta_\ell}(\xi)\overline{\widehat{\ta_\ell}(\xi+2\pi
\omega)}=\delta(\omega), \qquad \omega\in \dmfcg{\dm}.
\end{equation}
As already mentioned before, the above equations in
\eqref{cond:basic} only depend on $\dmfcg{\dm}$, which in turn only
depend on the lattice $\dm \dZ=\Gamma$. As showed in \cite[Corollary
3.4]{Han:laa}, there exist two integer matrices $E$ and $F$ such
that $|\det E|=|\det F|=1$ and
\[
\dm=E D F, \qquad D=\hbox{diag}(\df_1, \ldots, \df_m, 1, \ldots, 1),
\quad \df_1\ge \ldots \df_m>1.
\]
For the dilation matrix $\dn:=\hbox{diag}(\df_1, \ldots, \df_m)$, we
first construct a tensor product tight affine frame filter bank as
follows. For each $\df_n>1$, one can easily construct a
one-dimensional tight affine frame filter bank $\tu_\ell, \ell=1,
\ldots, \mpsi_{\df_n}$ such that
\begin{equation}\label{cond:basic:1d}
\sum_{\ell=1}^{\mpsi_{\df_n}}
\widehat{\tu_\ell}(\xi)\overline{\widehat{\tu_\ell}(\xi+2\pi
\omega)}=\delta(\omega), \qquad \omega\in \{0, \tfrac{1}{\df_n},
\ldots, \tfrac{\df_n-1}{\df_n}\}.
\end{equation}
For the construction of  (compactly supported) one-dimensional tight
affine frame filters, see
\cite{DHRS,Han:acha:98,Han:twf,Han:ACHA:09,RS97a}.

Now we construct a tensor product filter bank by
\[
U_{(\ell_1, \ldots, \ell_m)}(\beta_1, \ldots, \beta_d):=
u_{\ell_1}(\beta_1) \cdots u_{\ell_m}(\beta_m), \qquad \beta_1,
\ldots, \beta_d\in \Z.
\]

This generates a total of $\mpsi:=\mpsi_{\df_1} \cdot \ldots \cdot
\mpsi_{\df_m}$ filters. We reorder them as $U_1, \ldots, U_\mpsi$.
By \eqref{cond:basic:1d}, one can easily check that
\[
\sum_{\omega\in \dmfcg{\dm}}
\widehat{U_\ell}(\xi)\overline{\widehat{U_\ell}(\xi+2\pi
\omega)}=\delta(\omega), \qquad \xi\in \dR.
\]
Now we define
$\widehat{a_\ell}(\xi):=\widehat{U_\ell}(E^T\xi)$, $\xi\in \dR,
\ell=1, \ldots, \mpsi$.
Then we have
\begin{equation}\label{tp:filter}
\sum_{\omega\in \Omega_\dm}
\widehat{\ta_\ell}(\xi)\overline{\widehat{\ta_\ell}(\xi+2\pi
\omega)}=\delta(\omega), \qquad \xi\in \dR.
\end{equation}
In other words, using a tensor product construction, for any
sublattice of $\dZ$ generated by $\dm \dZ$, we can easily obtain a
set of (compactly supported) tight affine(-like) frame filters satisfying
\eqref{tp:filter}. It is very important to notice that the above
construction only depends on the lattice $\dm \dZ$ instead of the
dilation matrix $\dm$ itself. For more details on construction of
high-dimensional wavelet filters using the tensor product methods,
see \cite{Han:laa,Han:twf}.

\medskip

Next, we shall present a general construction of tight affine(-like) frame
filters to fulfill the need in this paper. Suppose that we are given
a group of dilation matrices $\dm_\ell, \ell=1, \ldots, r$. To
achieve directionality, as discussed above, only the dilation
matrices for the low-pass filters will be important. In other words,
if shear matrices are involved, i.e., direction-based matrices, they
play an important role for low-pass filter only and for high-pass
filters, the choice of the dilation matrices does not matter, since
no further decomposition will be performed for high-pass
coefficients.

We group these dilation matrices into subgroups according to their
lattices $\dm_\ell \dZ$: if the lattices $\dm_\ell=\dm_{\ell'}$ are
the same, then $\dm_\ell, \dm_{\ell'}$ are grouped into the same
group. For each dilation matrix in a subgroup, we only use a fixed
set of tight affine frame filters that are constructed for that
lattice. In other words, for a given lattice, we have a set of tight
affine frame filters and we have complete freedom in choosing the
direction-based matrices, for instance, shear matrices, to achieve
directionality as long as the resulting lattice is the same given
lattice.

Suppose now that we choose $N$ sets of such tight affine frame filters
for all the groups of dilation matrices. Since these $N$ sets of
tight affine frames are completely independent, when we put them
together to get one whole set of tight affine frame filters, we have
to multiply the factor $\frac{1}{\sqrt{N}}$ to every involved filter
in the set. Now it is straightforward to see that the total
collection of all such $N$ sets of renormalized tight affine frame
filters indeed forms a collection of tight affine frame filters
satisfying \eqref{filter:general}, where $r$ is the total number of
all the involved filters. It is also very important to notice that
the renormalization of the filters does not reduce the
directionality of the tight affine frame system, since the
coefficients in the same band has the same ordering of magnitude as
the one without renormalization.

\subsection{Shearlet Constructions}

We next focus on explicit shearlet constructions. For this, we
recall that the dilation for shearlets is composed of a shear matrix
$S_s$ and a parabolic scaling matrix $A_c$ (cf. Subsection
\ref{subsec:shearlets}). Thus, at each level of the decomposition,
the matrices $\dm_\mi$ will be chosen as \beq
\label{eq:choiceMshearlet} \dm_\mi = S_\mi A_4, \eeq for some shear
matrices $S_\mi$ whose selection should be driven by the particular
application at hand, and $A_4$ is the parabolic scaling matrix,
where the choice of the value $4$ (in contrast to $2$) avoids
technicalities caused by square roots. As already elaborated upon
before, which of those matrices will be labeled low- and which one
high-pass is left to the user.

\subsubsection{Shearlet Unitary Extension Principle and Fast Shearlet Transform}

The requirements on the filters associated with general matrices of
the type described in \eqref{eq:choiceMshearlet} at each level were
already derived in Theorem \ref{theo:generalshearletUEP}, which we
coined the `Shearlet Unitary Extension Principle'.

Aiming towards examples for possible filters, we first observe that
there exists only two different lattices in
Theorem~\ref{theo:generalshearletUEP}: Consider the product $S_\ell A_4$ with
$S_\ell=\left(\begin{matrix} 1 &s_\ell\\ 0 &1\end{matrix}\right)$.
If $s_\ell$ is an even integer, then $S_\ell A_4 \Z^2=A_4 \Z^2$, if
it is an odd integer, then $S_\ell A_4
\Z^2=SA_4 \Z^2$ with $S=\left(\begin{matrix} 1 &1\\ 0
&1\end{matrix}\right)$. Hence we only need to design two sets of
tight shearlet frame filters in advance for the lattices $A_4 \Z^2$
and $SA_4 \Z^2$, respectively, following the general construction in
Subsection \ref{subsec:generalconstr}. Then we have complete freedom
in choosing the shear matrices $S_\ell$ to obtain a whole set of
tight shearlet frame filters satisfying perfect reconstruction (see
condition (i) in Theorem~\ref{theo:generalshearletUEP}). In
particular, notice that we can choose the filters such that the
designed tight shearlet frame is compactly supported, for example,
take $\dm_0=2I_2$ and choose the tensor product 1D tight wavelet
frames derived from the linear hat function as constructed in
\cite{RS97a}.

\begin{figure}[ht]
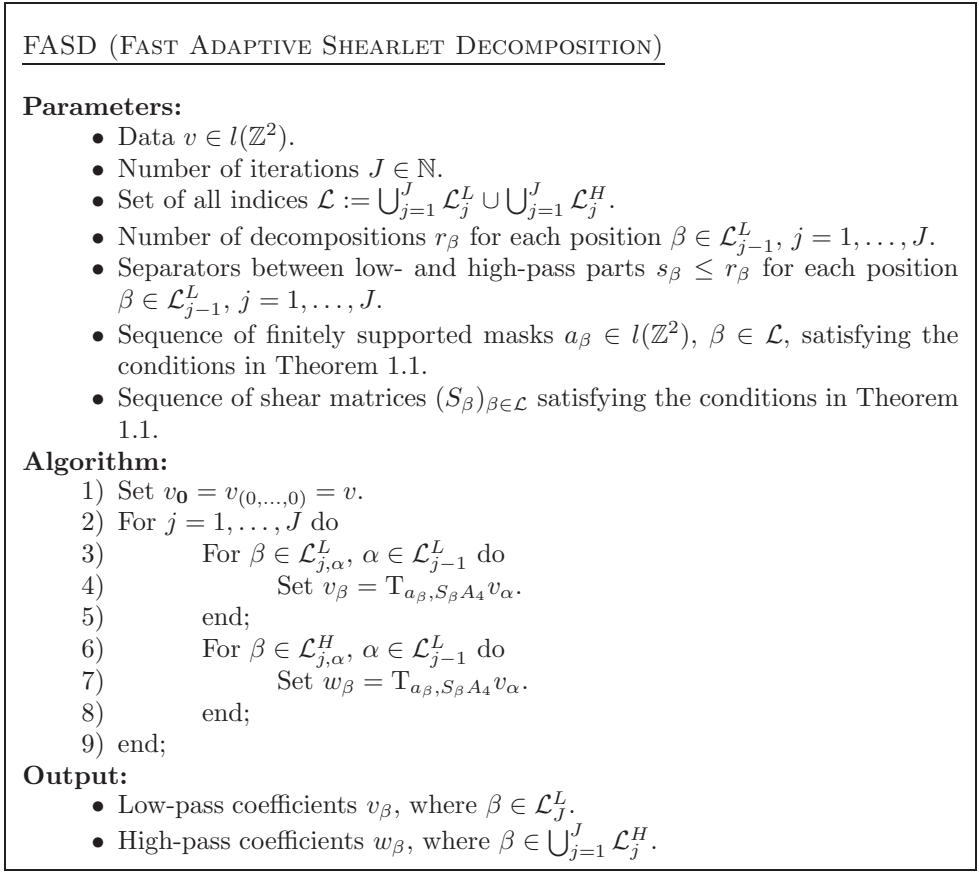

\centering \framebox{
\begin{minipage}[h]{4.9in}
\vspace*{0.3cm} {\sc \underline{FASD (Fast Adaptive Shearlet
Decomposition)}}

\vspace*{0.4cm}

{\bf Parameters:}\\[-3ex]
\begin{itemize}
\item Data $\tv \in \llll(\ZZ^2)$.
\item Number of iterations $J \in \NN$.
\item Set of all indices  $\cL:=\bigcup_{j=1}^J \cL^L_{j} \cup \bigcup_{j=1}^J \cL^H_{j}$.
\item Number of decompositions $\mpsi_\beta$ for each position $\beta \in \cL^L_{j-1}$, $j=1,\ldots,J$.
\item Separators between low- and high-pass parts $\mphi_\beta \le \mpsi_\beta$
for each position $\beta \in \cL^L_{j-1}$, $j=1,\ldots,J$.
\item Sequence of finitely supported masks $\ta_{\beta} \in \llll(\ZZ^2)$, $\beta \in \cL$,
satisfying the conditions in Theorem \ref{theo:generalshearletUEP}.
\item Sequence of shear matrices $(S_\beta)_{\beta \in \cL}$
satisfying the conditions in Theorem \ref{theo:generalshearletUEP}.
\end{itemize}

{\bf Algorithm:}\\[-3ex]
\begin{itemize}
\item[1)] Set $v_{\mathbf{0}} = \tv_{(0,...,0)} = \tv$.
\item[2)] For $j=1, \ldots, J$ do
\item[3)] \hspace*{1cm} For $\beta \in \cL^L_{j,\alpha}$, $\alpha \in \cL^L_{j-1}$ do
\item[4)] \hspace*{2cm} Set $\tv_\beta = \tz_{\ta_\beta, S_{\beta}A_4} \tv_{\alpha}$.
\item[5)] \hspace*{1cm} end;
\item[6)] \hspace*{1cm} For $\beta \in  \cL^H_{j,\alpha}$, $\alpha \in \cL^L_{j-1}$ do
\item[7)] \hspace*{2cm} Set $\tw_\beta = \tz_{\ta_\beta, S_{\beta}A_4} \tv_{\alpha}$.
\item[8)] \hspace*{1cm} end;
\item[9)] end;
\end{itemize}

{\bf Output:}\\[-3ex]
\begin{itemize}
\item Low-pass coefficients
$\tv_{\beta}$, where $\beta \in \cL^L_{J}$.
\item High-pass coefficients
$\tw_{\beta}$, where $\beta \in \bigcup_{j=1}^J \cL^H_{j}$.
\end{itemize}
\vspace*{0.01cm}
\end{minipage}
} \caption{The {\sc FASD} Algorithm for a fast adaptive
decomposition using shearlet systems.} \label{fig:fadshearlet}
\end{figure}

\begin{figure}[ht]
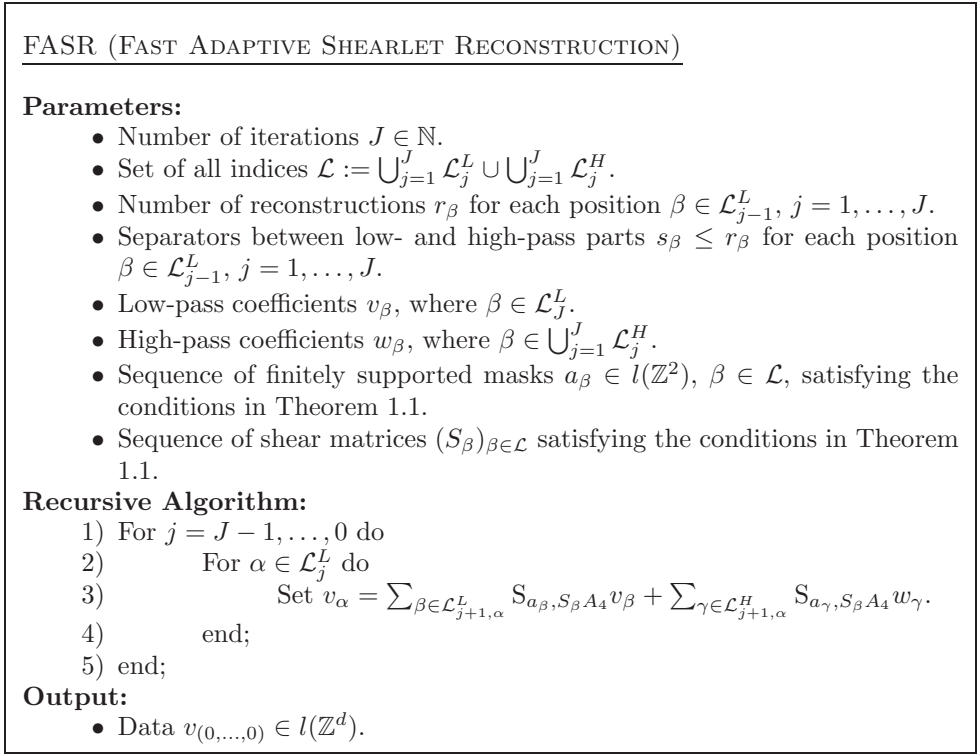

\centering \framebox{
\begin{minipage}[h]{4.9in}
\vspace*{0.3cm} {\sc \underline{FASR (Fast Adaptive Shearlet
Reconstruction)}}

\vspace*{0.4cm}

{\bf Parameters:}\\[-3ex]
\begin{itemize}
\item Number of iterations $J \in \NN$.
\item Set of all indices $\cL:=\bigcup_{j=1}^J \cL^L_{j} \cup \bigcup_{j=1}^J \cL^H_{j}$.
\item Number of reconstructions $\mpsi_\beta$ for each position $\beta \in \cL^L_{j-1}$, $j=1,\ldots,J$.
\item Separators between low- and high-pass parts $\mphi_\beta \le \mpsi_\beta$
for each position $\beta \in \cL^L_{j-1}$, $j=1,\ldots,J$.
\item Low-pass coefficients $\tv_{\beta}$, where $\beta \in \cL^L_{J}$.
\item High-pass coefficients $\tw_{\beta}$, where $\beta \in \bigcup_{j=1}^J \cL^H_{j}$.
\item Sequence of finitely supported masks $\ta_{\beta} \in \llll(\ZZ^2)$, $\beta \in \cL$,
satisfying the conditions in Theorem \ref{theo:generalshearletUEP}.
\item Sequence of shear matrices $(S_\beta)_{\beta \in \cL}$
satisfying the conditions in Theorem \ref{theo:generalshearletUEP}.
\end{itemize}

{\bf Recursive Algorithm:}\\[-3ex]
\begin{itemize}
\item[1)] For $j=J-1, \ldots, 0$ do
\item[2)] \hspace*{1cm} For $\alpha \in \cL^L_{j}$ do
\item[3)] \hspace*{2cm} Set $\tv_{\alpha} = \sum_{\beta \in \cL_{j+1,\alpha}^L} \sd_{\ta_{\beta}, S_{\beta}A_4} \tv_{\beta}
+\sum_{\gamma \in \cL_{j+1,\alpha}^H} \sd_{\ta_{\gamma},
S_{\beta}A_4} \tw_{\gamma}$.
\item[4)] \hspace*{1cm} end;
\item[5)] end;
\end{itemize}

{\bf Output:}\\[-3ex]
\begin{itemize}
\item Data $\tv_{(0,\ldots,0)} \in \llll(\ZZ^d)$.
\end{itemize}
\vspace*{0.01cm}
\end{minipage}
} \caption{The {\sc FASR} Algorithm for a fast adaptive
reconstruction using shearlet systems.} \label{fig:farshearlet}
\end{figure}

Each such a choice of filters for each level then leads to a Fast
Adaptive Shearlet Decomposition associated with a Fast Adaptive
Shearlet Reconstruction. Those two algorithms are described in
Figures \ref{fig:fadshearlet} and \ref{fig:farshearlet}. Notice that
here -- as already in the general version of {\sc (FAD)} and {\sc
(FAR)} -- the shear matrices can be chosen differently at each level
of the decomposition. We can envision that this adaption can be made
flexible dependent on a quick analysis, thresholding, say, of the
data outputted in the previous decomposition step. The algorithm
leaves all those possibilities open. The great flexibility provided
here should be utilizable for various applications including,
in particular, frame based image restorations.

\subsection{3D Shearlets}
\label{subsec:3D}

Various applications such as the problem of geometric separation in
image processing, for example, in astronomy when images of
galaxies require separated analyzes of stars, filaments, and sheets
call for directional representation systems for 3D data (see, e.g.,
\cite{DK08a,DK08aa}). In this subsection, we will show that our
framework conveys the potential to generate (compactly supported)
shearlets for 3D signals alongside a fast flexible decomposition
strategy.

In fact, our framework also provides shearlets associated with an
AMRA structure in higher dimensions. However, for the purpose of
clarity we will focus on dimension $3$. One possibility to extend
the parabolic scaling is by defining \beq \label{eq:para3D} A_c =
\left(\begin{array}{rrr} c & 0 & 0\\ 0 & \sqrt{c}& 0\\ 0 & 0 &
\sqrt{c}\end{array} \right),\quad c > 0. \eeq General shear matrices
were for instance studied in \cite{GLLWW06}. But in order to present
the idea for a 3D shearlet construction and don't cloud it with
technicalities, the shear matrices introduced in \cite{DST09} shall
serve as an interesting example, and we define \beq
\label{eq:shear3D} S_{s_1,s_2} = \left(\begin{array}{rrr} 1 & s_1 &
s_2\\ 0 & 1& 0\\ 0 & 0 & 1\end{array} \right),\quad (s_1,s_2) \in
\RR^2. \eeq We now set $ \dm_\mi = S_\mi A_4, $ for some shear
matrices $S_\mi$ of type \eqref{eq:shear3D}, whose selection should
be driven by the particular application at hand, and $A_4$ is the
parabolic scaling matrix from \eqref{eq:para3D}.

Concerning a construction of filters, we can again employ the
general construction from Subsection \ref{subsec:generalconstr}. By
calculation, it is easy to see that we have only four different sets
of sublattices, according to even or odd $s_1$ and $s_2$. Using
tensor products, we can design four sets of tight affine frame
filters in advances. Then we have the complete freedom in choosing
the shear matrices $S_{s_1, s_2}$ to obtain a whole set of tight
frame filters.

\medskip

We would like to mention that although focussing on 3D in this
subsection, our framework can also be used to generate shearlet
systems in higher dimensions using the general definition of shear
matrices in \cite{GLLWW06}. The parabolic scaling can be extended to
higher dimension, for instance, as it was done in \cite{DST09}, but
various other choices are possible, which will be explored in future
work.


\section{Further Examples of Directional Systems associated with an AMRA}

To show the broad applicability of the AMRA-structure, in this
section, we will present further examples of systems which have the
potential to analyze the directional content of images and signals.

\subsection{Shearlet-Like Construction}

The resolution of anisotropic structures requires paramatrization of
directions. The most natural way seems to be the choice of
parametrization through angles. However, angles are related with the
rotation operator, which is incompatible with digital lattices. This
observation led to the introduction of shearlets, which measure
directionality by slope rather than angle thereby greatly supporting
the treating of the digital setting. However, it seems natural to
seek for more general operators which leave lattice structures
invariant while being somewhat close to rotation. These might differ
from the shear matrices utilized in the shearlet definition, and
provide an improvement on uniform treatment of all angles. We would
like to remind the reader that in the shearlet construction (cf.
Subsection \ref{subsec:shearlets}) the partition into cones
serves as a means to allow a more uniform treatment of the angles,
since the `pure' shearlet system generated by the shearlet group is
strongly biased towards one axis (for more details we refer the
interested reader to \cite{GKL06}).

We can model this optimization problem in the following way:
\[
\mbox{{\sc (Rot$_\theta$)}} \quad B_\theta = \mbox{ argmin }_{B \in
M(\ZZ,2)} \int_0^{2\pi} \|(R_\theta - B) \left(\begin{array}{c} \cos
t\\ \sin t \end{array} \right)\|_2^2 \, dt \quad \mbox{s.t. }
\det(B)=1,
\]
where $R_\theta$, $\theta \in [0,2\pi)$ is planar rotation by
$\theta$ radians.

This optimization problem can be explicitly solved, which is done in
the following result:
\begin{theorem}
\label{theo:rotation} For $\theta \in [0,2\pi)$, we denote the
solution of {\sc (Rot$_\theta$)} by $B_\theta$. Then
\[
B_\theta = \left\{ \begin{array}{ccl}
I & : & 2\pi - \arcsin(\frac12) \le \theta < \arcsin(\frac12),\\
S_{-1}\mbox{ or }S_1^T
& : & \arcsin(\frac12) \le \theta < \frac{\pi}{4},\\
S_1^{(1)}\mbox{ or }S_1^{(2)}
& : & \frac{\pi}{4} < \theta \le \arccos(\frac12),\\
S_0^{(1)}
& : & \arccos(\frac12) < \theta \le \pi - \arccos(\frac12),\\
S_{-1}^{(1)}\mbox{ or }S_{-1}^{(2)}
& : & \pi - \arccos(\frac12) < \theta \le \frac{3\pi}{4},\\
-S_1\mbox{ or }-S_{-1}^T
& : & \frac{3\pi}{4} < \theta \le \pi - \arcsin(\frac12).
\end{array} \right.
\]
The remaining six cases can be determined by transposition of the
previous matrices, in particular,
\[
B_\theta = \left\{ \begin{array}{ccl}
-I & : & \pi - \arcsin(\frac12) < \theta \le \pi +
\arcsin(\frac12),\\
-S_1^T\mbox{ or }-S_{-1}
& : & \pi + \arcsin(\frac12) < \theta \le \frac{5\pi}{4},\\
(S_{-1}^{(1)})^T\mbox{ or }(S_{-1}^{(2)})^T
& : & \frac{5\pi}{4} < \theta \le \pi + \arccos(\frac12),\\
(S_0^{(1)})^T
& : & \pi + \arccos(\frac12) < \theta \le \frac{3\pi}{2} + \arcsin(\frac12),\\
(S_{1}^{(1)})^T \mbox{ or }(S_{1}^{(2)})^T
& : & \frac{3\pi}{2} + \arcsin(\frac12) < \theta \le \frac{7\pi}{4},\\
S_{-1}^T\mbox{ or }S_1
& : & \frac{7\pi}{4} \le \theta < 2 \pi - \arcsin(\frac12),\\
\end{array} \right.
\]
where
\[
S_s^{(1)}=\left(\begin{array}{rr} s & -1\\ 1 & 0 \end{array} \right)
\quad \mbox{and} \quad S_s^{(2)}=\left(\begin{array}{rr} 0 & -1\\ 1
& s \end{array} \right).
\]
\end{theorem}

\begin{proof}
Set $B = (b_{ij})_{1 \le i,j \le 2}$. First,
\begin{eqnarray*}
\lefteqn{\|(R_\theta - B) \left(\begin{array}{c} \cos t\\ \sin t \end{array} \right)\|_2^2}\\[1ex]
& = & [(\cos \theta - b_{11})^2 + (\sin \theta - b_{21})^2] \cos^2 t
+ [(-\sin \theta - b_{12})^2 + (\cos \theta - b_{22})^2] \sin^2 t\\
& & + 2 [(\cos \theta - b_{11})(- \sin \theta - b_{12}) + (\sin
\theta - b_{21})(\cos \theta - b_{22})]\cos t \,\sin t.
\end{eqnarray*}
Hence
\begin{eqnarray*}
\lefteqn{\int_0^{2\pi} \|(R_\theta - B) \left(\begin{array}{c} \cos t\\ \sin
t \end{array} \right)\|_2^2 \, dt}\\
& = & \pi [(\cos \theta - b_{11})^2 +
(-\sin \theta - b_{12})^2 + (\sin \theta - b_{21})^2 + (\cos \theta
- b_{22})^2].
\end{eqnarray*}
We now aim to minimize the term on the RHS over $b_{11}, b_{12},
b_{21}, b_{22} \in \ZZ$ under the constraint that $\det B = 1$.
Certainly, the solution to \beq \label{eq:reducedmin} \min_{B =
(b_{ij})_{i,j} \in M(\ZZ,2)} [(\cos \theta - b_{11})^2 + (-\sin
\theta - b_{12})^2 + (\sin \theta - b_{21})^2 + (\cos \theta -
b_{22})^2], \eeq by denoting the elements in $R_\theta$ by
$r^\theta_{ij}$, is
\[
b_{ij} = \rd(r^\theta_{ij}) \quad \mbox{for all } i,j \in \{1,2\}.
\]
where $\rd(\cdot)$ is the classical rounding to the nearest integer.
Notice that $b_{ij} \in \{-1,0,1\}$. The next closest solution
(possibly not unique) can be determined by choosing $i_0, j_0 \in
\{1,2\}$ such that
\[
\|\rd(r^\theta_{i_0,j_0}) - r^\theta_{i_0,j_0}\|_2 \ge
\|\rd(r^\theta_{i,j}) - r^\theta_{i,j}\|_2 \quad \mbox{for all } i,j
\in \{1,2\}.
\]
Then choose
\[
b_{i_0,j_0} = \left\{ \begin{array}{ccl} \sgn(b_{i_0,j_0}) \cdot 1 &
: & \rd(b_{i_0,j_0}) = 0,\\ 0 & : & \rd(b_{i_0,j_0}) \neq
0\end{array} \right.
\]
and
\[
b_{ij} = \rd(r^\theta_{i,j}) \quad \mbox{for all } (i,j) \neq
(i_0,j_0).
\]
In fact, we will see that provided the solution to
\eqref{eq:reducedmin} does not satisfies the constraint $\det B =
1$, then the just described closest one does, hence is a solution of
{\sc (Rot$_\theta$)}.

We now consider the case $0 \le \theta < \frac{\pi}{2}$. The other
cases can be handled similarly. \bitem
\item If $0 \le \theta < \arcsin^{-1}(\frac12)$, then the closest solution is $R_\theta$, which satisfies
the constraint.\\
\item If $\arcsin^{-1}(\frac12) \le \theta < \arccos^{-1}(\frac12)$, then the solution to \eqref{eq:reducedmin} is
\[
A=  \left(\begin{array}{rr} 1 & -1 \\ 1 & 1 \end{array} \right),
\]
for which $\det A = 2$. The closest one to this -- as just described
-- is
\[
B=\left(\begin{array}{rr} 1 & -1\\ 0 & 1 \end{array} \right) \mbox{ or
} \left(\begin{array}{rr} 1 & 0\\ 1 & 1 \end{array} \right) \quad
\mbox{for } \theta < \frac{\pi}{4},
\]
and
\[
B=\left(\begin{array}{rr} 0 & -1\\ 1 & 1 \end{array} \right) \mbox{ or
} \left(\begin{array}{rr} 1 & -1\\ 1 & 0 \end{array} \right) \quad
\mbox{for } \theta \ge \frac{\pi}{4}.
\]
These integer matrices obviously satisfy the constraint $\det B= 1$.\\
\item If $\arccos^{-1}(\frac12) \le \theta < \frac{\pi}{2}$, then the closest solution is $R_\theta$, which satisfies
the constraint. \eitem This proves the theorem.
\end{proof}

It is interesting to notice that all these matrices are in fact
shear matrices or closely related to shear matrices.

This consideration now enables us to propose a shearlet-like system
with presumable improved uniformity of the treatment of all angles.
For this, we let $\theta_\ell$, $\ell=1,\ldots,12$ be
representatives of the angles giving rise to the twelve different
cases discussed in Theorem \ref{theo:rotation}. We choose
\[
\dm_\mi = B_{\theta_\ell} A_4,\quad \ell=1,\ldots,12,
\]
for each level leaving the `cut' between low- and high-pass parts to
the user. The sequence of associated filters $a_\ell$,
$\ell=1,\ldots,12$ can be chosen similarly as in Section
\ref{sec:shearlets}.

\subsection{An Isotropic System for detecting Anisotropic Phenomena}

The following is some explanation for the tight frame in Definition \ref{defi:ASJ}
about how we can achieve directionality by using the simplest
system: dilation $\dm=2I_2$ and tensor product filters, but using
shear matrices in the middle. It is very important to notice that
for this case, the construction of all the filters in
Corollary~\ref{coro:generalUEPspecial} (ii) is independent of the
choice of the shear matrix. This simple fact has been noticed in
\cite{Han:laa} and makes the construction of multivariate wavelets
very simple.

From now on, we assume that $\dm$ is fixed (say, $\dm=2I_d$) and all
filters are obtained by a tensor product construction. Therefore,
these filters give preferences to horizontal and vertical direction.
In other words, they are suitable for handling only edges having
either horizontal or vertical orientation.

Let us look at the elements in $\cAS_J$ introduced in Definition \ref{defi:ASJ}. For the
functions in level $1$, we have $\alpha\in \cL^L_1\cup \cL^H_1$,
that is, $\alpha=(\mi,0, \ldots,0)$ with $0<\mi\le \mpsi_\alpha$.
Then $\psi^\alpha_{\dn_\beta \dm_0^J; \vn}$ are given by
\begin{eqnarray*}
\psi^\alpha(\dm_\alpha^{-1} \dm_0^J\cdot-\vn)
& = & |\det (\dm_\alpha)| \sum_{\vk\in \dZ} \ta_\alpha(\vk)
\phi(\dm_0^J\cdot-\dm_\alpha\vn-\vk)\\
& =& |\det (\dm_0)| (\ta_\alpha
\star \phi)(\dm_0^J \cdot-\dm_\alpha\vn).
\end{eqnarray*}
Note that $\psi^\beta_{\dn_\beta \dm_0^J; \vn}$ with a high-pass
filter has the same directionality as the filter $\ta_\beta$, that
is, only either horizontal or vertical direction. Note that we
assume to have horizontal or vertical directionality for
high-pass filters only, while for low-pass filters, we assume that
they are isotropic and there is no directionality for them. In other
words, for low-pass filters, the filters have no directionality.

The directionality comes at the next level starting with level $2$.
For $\beta\in \cL^L_2\in \cL^H_2$, we have
$\beta=\alpha+(0,\mi,0,\ldots, 0)$ with $\alpha\in \cL^L_1$ and
$0<\mi\le \mpsi_\beta$. Note $\dn=\dm_\beta^{-1} \dm_{\alpha}^{-1}$.
So, $\psi^\beta_{\dn_\beta \dm_0^J; \vn}$ are given by
\[
\psi^\beta_{\dn_\beta \dm_0^J; \vn}= \psi^\beta(\dm_\beta^{-1}
\dm_{\alpha}^{-1} \dm_0^J\cdot-\vn) =|\det (\dm_0)|\sum_{\vk\in \dZ}
\ta_\beta(\vk) \psi^\alpha(\dm_\alpha^{-1}
\dm_0^J\cdot-\dm_\beta\vn-\vk).
\]
That is, we similarly have
\[
\psi^\beta_{\dn_\beta \dm_0^J; \vn}=|\det (\dm_0)| (\ta_\beta\star
\psi^\alpha)(\dm_\alpha^{-1} \dm_0^J\cdot-\dm_\beta\vn).
\]

To see the directionality better, we rewrite it in terms of $\phi$
as follows:
\begin{align*}
\psi^\beta_{\dn_\beta \dm_0^J; \vn}
&=|\det (\dm_0)|\sum_{\vk\in \dZ} \ta_\beta(\vk) \psi^\alpha(\dm_\alpha^{-1} \dm_0^J\cdot-\dm_\beta\vn-\vk)\\
&=|\det (\dm_0)|^2 \sum_{\vk\in \dZ} \ta_\beta(\vk) (\ta_\alpha\star
\phi)(\dm_0^J\cdot-\dm_\alpha \dm_\beta\vn -\dm_\alpha \vk).
\end{align*}
To see the directionality of $\psi^\beta_{\dn_\beta \dm_0^J; \vn}$,
it suffices to see it for $\vn=0$, since all others are just shifts
of $\vn=0$ (that is, they have the same directionality). In this
case,
\[
\psi^\beta_{\dn_\beta \dm_0^J; 0}=|\det (\dm_0)|^2 \sum_{\vk\in \dZ}
\ta_\beta(\vk) (\ta_\alpha\star \phi)(\dm_0^J\cdot -\dm_\alpha \vk).
\]
Note that $\ta_\alpha$ is a low-pass filter and therefore, it has no
directionality. In other words, the function $\ta_\alpha\star \phi$
has most of its energy concentrated about the
origin. Though $\dm_0^J$ can also yield directionality by rotating
the function, let us ignore this possibility here and simply assume
that $\dm_0=2I_d$ (or even simply take $J=0$).

Now let us argue that if $\ta_\beta$, as a high-pass filter, has
directionality either in the horizontal or vertical direction, then
the function $\psi^\beta_{\dn_\beta \dm_0^J; 0}$ will have
directionality other than horizontal and vertical, induced by
$\dm_\alpha= S_\alpha \dm_0$, where $S_\alpha$ is some shear matrix.
For simplicity, we assume $\ta_\beta$ is a filter with horizontal
directionality. Therefore, $\ta_\beta(\vk)$ concentrates most of its
energy for $\vk$ on/around the $x$-axis. But the energy of the term
$(\ta_\alpha\star \phi)(\dm_0^J\cdot -\dm_\alpha \vk)$ concentrates
around the point $\dm_0^{-J} \dm_\alpha \vk$. Therefore, for
$\vk=(k, 0, \ldots, 0)$ with $k\in \ZZ$ (this is the line where most
energy of $\ta_\beta$ is concentrating), the energy location for
$\psi^\beta_{\dn_\beta\dm_0^J; 0}$ in above sum is put at the
location $\dm_0^{-J} \dm_\alpha \vk=\dm_0^{1-J} S_\alpha \vk$
instead of $\dm_0^{1-J} \vk$. In other words, the energy of the
function $\psi^\beta_{\dn_\beta\dm_0^J; 0}$ is aligned around the
straight line: the image of the $x$-axis under the mapping
$\dm_0^{1-J} S_\alpha$. Since $S_\alpha$ acts like a rotation, the
directionality is controlled by this mapping; the
horizontal and vertical lines are mapped into other directions. By
using different $S_\alpha$ (such as the power of some shear matrices),
as many directions as required by the application at hand can be
created.


\end{document}